\documentclass[11pt]{amsart}
\usepackage{amsmath}
\usepackage{amsthm}
\usepackage{amsfonts}
\usepackage{amssymb}
\usepackage{epsfig}
\usepackage{color}
\usepackage{tikz}
\usepackage{enumitem} 
\usepackage{pgfplots}
\pgfplotsset{compat=newest}
\usetikzlibrary{plotmarks}
\usepackage{xcolor,colortbl}
\usepackage{multirow}
\usepackage{rotating}
\usepackage{booktabs}
\usepackage{appendix}
\usepackage{bm,bbm}
\usepackage{xfrac}
\usepackage[hidelinks]{hyperref}

\usepackage{mathtools}  
\usepackage[T1]{fontenc}
\newtheorem{theorem}{Theorem}
\newtheorem{proposition}[theorem]{Proposition}

\theoremstyle{definition}

\newtheorem{example}[theorem]{Example}
\theoremstyle{lemma}
\newtheorem{lemma}[theorem]{Lemma}

\theoremstyle{remark}
\newtheorem{remark}[theorem]{Remark}

\newtheorem{assumption}[theorem]{Assumption}
\numberwithin{theorem}{section}
\numberwithin{equation}{section}
\numberwithin{table}{section}
\numberwithin{figure}{section}
\input{def}
\allowdisplaybreaks
\begin{document}
\title[Singular perturbation results for linear PDAEs of hyperbolic type]{Singular perturbation results for linear partial differential--algebraic equations of hyperbolic type$^*$}
\author[]{R.~Altmann$^\dagger$, C.~Zimmer$^{\dagger}$}
\address{${}^{\dagger}$ Department of Mathematics, University of Augsburg, Universit\"atsstr.~14, 86159 Augsburg, Germany}
\email{robert.altmann@math.uni-augsburg.de, christoph.zimmer@math.uni-augsburg.de}
\thanks{${}^{*}$ Research funded by the Deutsche Forschungsgemeinschaft (DFG, German Research Foundation) -- Project-ID 446856041. 	 
}
\date{\today}
\begin{abstract}
We consider constrained partial differential equations of hyperbolic type with a small parameter $\eps>0$, which turn parabolic in the limit case, i.e., for $\eps=0$. The well-posedness of the resulting systems is discussed and the corresponding solutions are compared in terms of the parameter~$\eps$. For the analysis, we consider the system equations as partial differential--algebraic equation based on the variational formulation of the problem. For a particular choice of the initial data, we reach first- and second-order estimates. For general initial data, lower-order estimates are proven and their optimality is shown numerically. 
\end{abstract}
%
%
\maketitle
%
%
%
\section{Introduction}\label{sect:introduction}
Singularly perturbed differential and evolution equations have been analyzed for many decades already (see~\cite{KadP03} for a review) and cover the entire spectrum of elliptic, parabolic, 
as well as hyperbolic systems. In this paper, we are concerned with the asymptotics of hyperbolic--parabolic singular systems, i.e., hyperbolic systems with a parabolic limit case. These systems are well-studied in the unconstrained case: In~\cite{Lio73,GhiG06}, linear and rooty convergence rates for the variables of the linear damped wave equation are proven. 
Later, these results were extended to nonlinear systems \cite{EshW88,Esh90,HasY07,GhiG12} as well as integral expressions~\cite{Esh87b,LiaLX05}. 
Resulting estimates can be used, e.g., in the theory of boundary layers or as a tool for the design of numerical algorithms~\cite{Esh87}, which also serves as motivation for the present work.

For the particular case of one-dimensional gas networks, i.e., a coupled system of hyperbolic equations, perturbation results have been derived in~\cite{EggK18}. In general, such network structures can be modeled as constrained partial differential equations, where the constraints are naturally given by the junctions within the network, reflecting fundamental physical properties. Such an approach leads to {\em partial differential--algebraic equations} (PDAEs), cf.~\cite{EmmM13,LamMT13,Alt15}, which may be interpreted as differential--algebraic equations in Banach spaces. For the sake of completeness, we would like to mention that one may also consider the network as a domain on which the differential equations are stated~\cite{Mug14}. In this approach, however, inhomogeneous boundary conditions still account for constraints on the solution. 

In this paper, we focus on singularly perturbed (linear) PDAEs of first order (in time) that are hyperbolic, i.e., we consider hyperbolic partial differential equations including a small parameter~$\eps$, which underlie an additional constraint. 
The singular perturbation of the considered PDAEs is characterized through the property that the system is of hyperbolic nature for~$\eps>0$ and parabolic in the limit case~$\eps=0$.
Throughout the paper, we present three examples in more detail. 
\begin{example}[Damped acoustic wave equation]\label{exp:wave}
Consider the damped wave equation~$\ddot u + d\dot u - c^2\Delta u = 0$ with the damping parameter~$d$ being proportional to the square of the speed of sound~$c$.  
Setting $\eps = 1/c^2$, we can rewrite this as 
\begin{align*}
  \dot u + \nabla\cdot w &= 0,\\
  \eps\, \dot w + \nabla u + \tfrac{d}{c^2} w &= 0,  
\end{align*}
cf.~\cite[Ch.~I.1]{Bra07}. In particular, $\eps$ is small (i.e., $c$ is large) if we consider the propagation of acoustic waves in solids~\cite{CRC13}. Note that, for the limit case $\eps=0$, this system reduces to the (parabolic) heat equation in mixed form with thermal diffusivity~$c^2/d$. 
In both cases, Dirichlet boundary conditions can be incorporated with the help of a Lagrange multiplier leading to a PDAE~\cite{HinPUU09,Alt15}.
\end{example}	
\begin{example}[Viscoelastic Stokes problem]\label{exp:Stokes}
Flows of fluids with complex microstructure, e.g., molten polymers, can be modeled by
\begin{align*}
  \rho\,\dot v - \nabla\cdot T + \nabla p &= \g,\\
  \eps\, \dot T - \eta\, \big(\nabla v + (\nabla v)^\top \big)+ T &= \f, 
\end{align*}
together with the constraint $\nabla \cdot v = 0$, cf.~\cite{Ren89}. Here, $v$ denotes the velocity field, for which we assume homogeneous Dirichlet boundary conditions, $T$ the Cauchy stress tensor, and~$p$ the pressure. Moreover, the density~$\rho$, the zero-shear-rate viscosity $\eta$, and the relaxation time~$\eps$ are positive constants. This model can be seen as the linearized Navier--Stokes equation for viscoelastic fluids. In particular, the system reduces to the unsteady Stokes equation for a vanishing relaxation time~$\eps$.
\end{example}
\begin{example}[Electro-magnetic energy propagation in power networks]\label{exp:pressure}
The third example describes the electro-magnetic energy propagation in power networks~\cite{MagWTA00,GoeHS16} as well as the propagation of pressure waves in a network of gas pipes~\cite{Osi87,BroGH11,JanT14,EggKLMM18}. 
The geometry of the underlying network can be encoded by a directed graph. In power networks, the edges represent transmission lines, whereas the vertices model a customer, a power supplier, or an interconnection. 

From an analytical point of view, it is sufficient to consider a single transmission line~\cite{EggK18}. Hence, we consider the unit interval as physical domain. Under certain simplifying model assumptions, the associated system in its strong form is given by
\begin{align*}
	\dot p + a\, p + \partial_x m &= \g,  \\ 
	\eps\, \dot m + \partial_x p + d\,m &= \f, 
\end{align*}
where $a$ and $d$ are (space-dependent) damping parameters describing the resistance and conductance in the transmission line. This system is also known as {\em telegrapher's equation}~\cite{GoeHS16}. 
Similar to the first example, constraints may occur due to inhomogeneous Dirichlet boundary conditions for the voltage~$p$ or as modeling tool to include the network structure. 

For the propagation of pressure waves, where we have $a=0$ (cf.~\cite{EggK17ppt,AltZ18ppt}), $\eps$ equals the product of the adiabatic coefficient and the square of the Mach number and is of order~$10^{-3}$, cf.~\cite{BroGH11}. 
\end{example}
The first main result of this paper compares the two corresponding solutions and shows that they only differ by a term of order~$\eps$, as long as the initial data is chosen appropriately. At this point, we would like to emphasize that such a condition does not occur in the respective finite-dimensional setting, which was analyzed in~\cite[Ch.~2.5]{KokKO99}. 
Also in the infinite-dimensional setting, one can renounce such a condition if the initial values are sufficiently regular~\cite[Ch.~2]{GhiG06}.  
For general initial data, on the other hand, we loose half an order, leading to an order of~$\sqrt{\eps}$ only.  
The second main contribution considers a second-order approximation of the original solution. Again, sufficient conditions on the initial data and regularity assumptions on the right-hand sides are discussed, which guarantee the full approximation order of two. 
Similarly as before, general initial data reduces the order to~$\eps^{\sfrac 32}$ or even~$\sqrt{\eps}$. To show that the presented estimates are indeed sharp, we examine specific numerical examples. 
\smallskip

The remaining parts of the paper are organized as follows. In Section~\ref{sect:PDAE} we introduce the functional analytic setting for linear first-order PDAEs of hyperbolic type, including a small parameter~$0<\eps\ll 1$. 
Moreover, we show that the particular examples introduced above fit into the presented framework. 
The existence of mild and classical solutions as well as stability estimates are then discussed in Section~\ref{sect:existence}. This also covers the discussion of the limit case for~$\eps=0$, which is of parabolic nature. 
The main results of the paper in Section~\ref{sect:expansion} are devoted to the comparison of the solutions of the original and the limit equations, leading to first and second-order estimates in terms of the parameter~$\eps$. Finally, the theoretical approximation orders a numerically verified by a number of experiments in Section~\ref{sect:numerics}.
\smallskip

Throughout this paper, we use for estimates the notion $a \lesssim b$ for the existence of a generic constant $c>0$ such that $a \le cb$.  
%
%
\section{Hyperbolic PDAE Model}\label{sect:PDAE}
In this section, we introduce the system class of interest namely first-order hyperbolic systems with a small parameter~$\eps$, which satisfy an additional constraint. 
Moreover, we discuss the already mentioned examples and show that they fit in the given framework.  
%
%
\subsection{Function spaces and system equations}\label{sect:PDAE:prelim} 
For a general formulation of constrained hyperbolic systems of first order, we introduce the three Hilbert spaces~$\P$, $\M$, and $\Q$. The corresponding solution components will be denoted by~$p$, $m$, and $\lambda$, respectively. 
In general, one may think of~$p$ modeling a {\em potential} whereas $m$ is a {\em flow variable}. The variable~$\lambda$ serves as Lagrange multiplier for the incorporation of the given constraint. 
We assume that $\P$ forms a Gelfand triple with pivot space~$\calH$, i.e., $\P\, {\hookrightarrow}\, \calH \cong \calH^*\, {\hookrightarrow}\, \P^*$ where all embeddings are dense. 
On the other hand, $\M$ is assumed to be identifiable with its own dual space, i.e., $\M\cong\M^*$. For the applications in mind, $\calM$ equals an~$L^2$-space, see Section~\ref{sect:PDAE:example}. 

As we consider time-dependent problems, appropriate solution spaces are given by Sobolev--Bochner spaces; see~\cite[Ch.~7]{Rou05} for an introduction. Denoting the space of quadratic Bochner integrable functions with values taken in a Banach space~$X$ by~$L^2(0,T;X)$, we use the notion~$H^m(0,T;X)$, $m\in\N$, for functions with higher regularity in time. Moreover, we define for two Sobolev spaces $X_1 \hookrightarrow X_2$ the space 
\[
  W(0,T;X_1,X_2) 
  := \big\{ v\in L^2(0,T;X_1)\ |\ \dot v \text{ exists in } L^2(0,T;X_2) \big\}. 
\]
Within this paper, we consider PDAEs of the form 
\begin{subequations}
\label{eq:PDAE:eps}
\begin{alignat}{5}
	\dot p&\ +\ &\A p &\ -\ &\K^* m & + \B^*\lambda\ &=&\ \g &&\qquad \text{in }\P^*, \label{eq:PDAE:eps:a}\\
	\eps\, \dot m&\ +\ &\K p &\ +\ &\D m & &=&\ \f  &&\qquad \text{in }\M^*, \label{eq:PDAE:eps:b}\\
	& & \B p & & & &=&\ \h &&\qquad \text{in } \Q^* \label{eq:PDAE:eps:c}
\end{alignat}
\end{subequations}
with initial conditions for~$p(0)$ and $m(0)$. 
Note that~\eqref{eq:PDAE:eps:a} and~\eqref{eq:PDAE:eps:b} are differential equations, whereas~\eqref{eq:PDAE:eps:c} reflects a constraint on~$p$, which defines the PDAE structure. Further note that all three equations are formulated in the respective dual spaces, which display the corresponding space of test functions. 
Within this setting, the operators~$\A$ and~$\D$ introduce damping, $\K$ is a differential operator, and $\B$ incorporates the constraint. 
The precise assumptions are summarized in the following and will be validated for the given examples in Section~\ref{sect:PDAE:example} below. 
\begin{assumption}
\label{ass:operators}	
All operators 
\[
  \A\colon \calH \to \calH^*, \qquad 
  \K\colon \P \to \M^*, \qquad
  \D\colon \M \to \M^*, \qquad 
  \B \colon \P\to \Q^*.
\]
are linear and bounded. The operator norm of $\A$ is denoted by $C_{\A} \geq 0$ and analogously for $\K$, $\D$, and $\B$. Further, we assume that the operator~$\D$ is elliptic with 
\[
  c_{\D}\, \| m\|_\M^2
  \le \langle \D m, m\rangle
\]
and that $\B$ is inf-sup stable, i.e., there exists a positive constant~$\beta$ with 
\[
   \adjustlimits\inf_{\mu\in\Q \setminus \{ 0\}} \sup_{q\in \P\setminus \{ 0\}} \frac{\langle \B q, \mu \rangle}{\|q\|_\P \|\mu\|_\Q} = \beta >0.
\]
Finally, a constant~$c_{\K}>0$ exists, such that for all $q_{\ker} \in \Pker := \ker \B \subseteq \P$ we have 
\begin{equation}
\label{eq:ass_K}
c_{\K}\,\|q_{\ker}\|_{\P} \le \|\K q_{\ker} \|_{\M^*}. 
\end{equation}
\end{assumption}
As usual for constrained systems, the inf-sup stability of the constraint operator~$\B$ is a crucial property for the well-posedness of the PDAE~\eqref{eq:PDAE:eps}. This is due to the saddle point structure of the system equations. Further note that~$\B$ is automatically inf-sup stable if it is surjective and $\Q^*$ is a finite-dimensional space; see~\cite{AltZ18ppt}. Contrariwise, the inf-sup condition implies surjectivity of the operator~$\B$ as well as injectivity of its dual~$\B^*$. The classical result presented in~\cite[Lem.~III.4.2]{Bra07} implies the existence of a right inverse~$\B^-\colon \Q^*\to\P$ and the estimates 
\[
  \|\B q\|_{\Q^*} \ge \beta\, \|q\|_\P 
  \qquad\text{and}\qquad 
  \|\B^- \h\|_\P \le \beta^{-1} \|\h\|_{\Q^*}
\]
for all~$q\in\Pker^{\perp}$ and $\h\in\Q^*$.  

The right-hand sides in~\eqref{eq:PDAE:eps} are of the form 
\[
\g\colon[0,T] \to \P^*, \qquad 
\f\colon[0,T] \to \M^*, \qquad 
\h\colon[0,T] \to \Q^*.
\]
Finally, the parameter~$\eps>0$ is expected to be small such that the term~$\eps\, \dot m$ in~\eqref{eq:PDAE:eps:b} takes the role of a singular perturbation. The limit case for~$\eps=0$ will be subject of Subsection~\ref{sect:existence:noEps}.  
An initial value of~$p$ is called {\em consistent} if the difference $p(0)-\B^- \h(0)$ is an element of the closure of $\Pker$ in $\calH$, which we denote by~$\cHker$ in the sequel. 
For~$p(0)\in\P$, the consistency conditions turns into~$\B p(0)= \h(0)$. 
%
%
\subsection{Examples}\label{sect:PDAE:example}
In this subsection we explain the meaning of the operators and spaces in the context of the three examples introduced in Section~\ref{sect:introduction}.
\subsection*{Example 1: Damped acoustic wave equation}
Recall the wave equation from Example~\ref{exp:wave} on a bounded domain~$\Omega \subset \R^n$ with Dirichlet boundary conditions for the state~$u$. 
For the corresponding PDAE formulation, we consider the spaces  
\[
  \P = H^1(\Omega), \qquad
  \cH = L^2(\Omega), \qquad
  \M = [L^2(\Omega)]^n, \qquad
  \Q = H^{-\sfrac 1 2}(\partial\Omega), 
\]
where the state $p$ denotes the density (called $u$ in the introduction) 
and $m$ the velocity. 
The operators are given by~$\A=0$, $\B$ denotes the trace operator to include Dirichlet boundary conditions, $\D$ equals the multiplication by the positive constant~$d/c^2$, and $\K$ is the weak gradient, meaning that $\langle\K q, m\rangle_{\calM^\ast,\calM} = \int_\Omega \nabla q \cdot m\dx$. In particular, the operator is $\B$ is inf-sup stable and the assumption~\eqref{eq:ass_K} on~$\K$ is satisfied for all $q_{\ker} \in \Pker = H^1_0(\Omega)$ by the Poincar\'e inequality \cite[Ch.~II.1]{Bra07}. Furthermore, it is well-known that $H^1_0(\Omega)$ is dense in $L^2(\Omega)$, such that $\cHker=\cH$. The right-hand sides are given by $\g=0$, $\f=0$, and $\h$ including the Dirichlet data. 
\subsection*{Example 2: Viscoelastic Stokes problem}
We consider the equations from Example~\ref{exp:Stokes} with a rescaling such that $\rho=1$ and $\eta=1/2$. Given a bounded domain~$\Omega \subset \R^n$, we define 
\[
  \P = [H^1_0(\Omega)]^n, \qquad
  \cH = [L^2(\Omega)]^n, \qquad
  \M = [L^2(\Omega)]_\text{sym}^{n,n}, \qquad
  \Q = L^2(\Omega) \setminus \R, 
\]
where $[L^2(\Omega)]_\text{sym}^{n,n}$ denotes the space of symmetric $n\times n$ matrices with coefficients in~$L^2(\Omega)$ and $L^2(\Omega) \setminus \R$ the space of $L^2$ functions with vanishing mean value. 
Here, $p$, $m$, and $\lambda$ denote the velocity field, the Cauchy stress tensor, and the pressure, respectively. 
For the operators, we define~$\A=0$ and~$\D$ as the identity. Moreover, $\B$ equals the (inf-sup stable) divergence operator such that~$\Pker$ is the space of divergence-free functions and~$\cHker$ its closure in $L^2$. For more details on these spaces, we refer to~\cite{Tar06}. Finally, $\K$ denotes the negative (weak) symmetric gradient, i.e., $\langle\K q, m\rangle_{\calM^\ast,\calM} = -\frac 12 \int_\Omega (\nabla q + (\nabla q)^\top ) : m\dx$, using the classical double dot notation from continuum mechanics.
\subsection*{Example 3: Electro-magnetic energy propagation in power networks} 
In this final example, we recall the equations of the propagation of the electro-magnetic energy in a transmission line, cf.~Example~\ref{exp:pressure}. The unknowns~$p$ and~$m$ equal the voltage and the current, respectively. For the weak formulation, we define 
\[
\P = H^1(0,1), \qquad
\M = \cH = L^2(0,1), \qquad
\Q = \R^2.  
\]
If we incorporate the boundary conditions for the potential~$p$ in form of a constraint, then~$\B$ equals the trace operator (point evaluation at the end points) and, hence, $\Pker=H^1_0(0,1)$ and $\cHker = \cH$, since $H^1_0(0,1)$ is dense in~$L^2(0,1)$. 
The operators $\A$ and $\D$ denote the multiplication by constants and include possibly state-dependent damping and friction to the model, whereas~$\K\colon \P\to\M^*$ equals the spatial derivative. With the help of the right-hand sides one can incorporate, e.g., the slope of a pipe. 

If the full power network is considered, then one may define $\P$ as the space of globally continuous and piecewise $H^1$-functions and $\M$ as the space of piecewise $L^2$-functions. In this setting, typical coupling conditions resemble the Kirchhoff circuit laws, which can be formulated as a right-hand side $\g$. For more details, we refer to~\cite{AltZ18ppt}. 
An alternative approach to the network case is considered in~\cite{Mug14}. 
%
%
\section{Existence of solutions}\label{sect:existence}
In order to prove the existence of mild and classical solutions as well as weak solutions for the parabolic limit case, we first discuss the solvability of a related stationary problem. 
%
%
\subsection{An auxiliary problem}\label{sect:existence:stationary}
It turns out that the following auxiliary problem is helpful for the upcoming analysis, 
\begin{subequations}
\label{eq:op:stationary}
\begin{alignat}{5}
	(\A+C_{\A} \id) \pb &\ -\ &\K^* \mb & + \B^*\lb &\ =\ &\g  &&\qquad \text{in }\P^*, \label{eq:op:stationary:a}\\
	\K \pb &\ +\ &\D \mb &   &\ =\ &  0 &&\qquad \text{in }\M^*, \label{eq:op:stationary:b}\\
	\B \pb & & &  &\ =\ &  0 &&\qquad \text{in }\Q^*. \label{eq:op:stationary:c}
\end{alignat}
\end{subequations}
Note that the system does not include time derivatives of the variables but that the right-hand side may still be time-dependent. To show the existence of a unique solution $(\pb,\mb,\lb)$ we first define an elliptic operator~$\wlap$.
\begin{lemma}
\label{lem_wlap}
Given Assumption~\ref{ass:operators}, the operator $\wlap:= \K^* \D^{-1} \K\colon \P \to \P^*$ is linear, continuous, and non-negative. 
Furthermore, its restriction to~$\Pker$ is elliptic, i.e., there exists a constant $\cLap>0$ such that for all $q_{\ker} \in \Pker$  we have that
\[ 
	\langle \wlap q_{\ker} , q_{\ker} \rangle \ge \cLap\, \| q_{\ker}\|^2_{\P}. 
\]
\end{lemma}
\begin{proof}
By the ellipticity of $\D$, the operator~$\wlap$ is well-defined and its linearity and continuity are obvious. 
For the non-negativity, we apply the ellipticity of~$\D$, leading to 
\begin{align*}
	\langle \wlap q, q \rangle_{\P^*,\P}  
	= \langle  \D^{-1}\K q, \K q \rangle_{\M,\M^*} 
	= \langle  \D^{-1}\K q, \D \D^{-1} \K q \rangle_{\M,\M^*} 
	\geq c_{\D} \| \D^{-1}\K q \|^2_{\M} 
	\geq 0 
\end{align*}
for all~$q\in \P$. In addition, for $q_{\ker}\in \Pker$ we can apply the properties of~$\K$. 
This yields
\begin{align*}
	\langle \wlap q_{\ker}, q_{\ker} \rangle_{\P^*,\P}  
	\geq c_{\D} \| \D^{-1}\K q_{\ker} \|^2_{\M} 
	\geq \frac{c_{\D}}{C_{\D}^2} \| \D \D^{-1} \K q_{\ker} \|^2_{\M^*} 
	\geq \frac{c_{\D}c_\K}{C_{\D}^2} \| q_{\ker} \|^2_{\P}, 
\end{align*}
which shows the claimed ellipticity. 
\end{proof}
Based on the newly introduced operator from the previous lemma, we now prove the existence of a solution to~\eqref{eq:op:stationary}. 
\begin{lemma}[Existence result for the auxiliary problem]
\label{lem_op_B_station}
Given Assumption~\ref{ass:operators} and a right-hand side $\g \in H^m(0,T;\P^*)$ for some $m\in\N$,  system~\eqref{eq:op:stationary} has a unique solution 
\[ 
  \big(\pb,\mb,\lb \big) 
  \ \in\  H^m(0,T; \Pker) \times H^m(0,T; \M) \times H^m(0,T;\Q),
\]
which depends continuously on the right-hand side~$\g$. 
\end{lemma}
\begin{proof}
We show that for~$\g \in \P^*$ (independent of time) system~\eqref{eq:op:stationary} has a unique solution
\[
  \big(\pb,\mb,\lb \big) \in \Pker \times \M \times \Q.   
\]
The result for a time-dependent right-hand side~$\g \in H^m(0,T;\P^*)$ then follows immediately by considering system~\eqref{eq:op:stationary} pointwise in time. 
The resulting solution is~$H^m$-regular in time, since all involved operators are time-independent.

Now consider $\g \in \P^*$. Since the operator~$\D$ is invertible, we can insert equation~\eqref{eq:op:stationary:b} into~\eqref{eq:op:stationary:a}, which results in the system  
\begin{alignat*}{5}
	(\wlap + \A+C_{\A} \id)\pb&\ + \B^*\lb &\ =\ &\g  &&\qquad \text{in } \P^*,\\
	\B \pb&  &\ =\ & 0 &&\qquad \text{in } \Q^*.
\end{alignat*}
By standard arguments~\cite[Ch.~II.1.1]{BreF91} this system has a unique solution $\pb \in \Pker$, $\lb \in \Q$, which is bounded in terms of $\g$.  
The existence of~$\mb$ and the stability bound~$ \|\pb\|_{\P} + \|\mb\|_\M + \| \lb \|_\Q \lesssim \|\g\|_{\P^*}$ then follow by $\mb = - \D^{-1} \K \pb$. 
\end{proof}
As another preparation for the existence results in the upcoming subsection, we consider the following lemma. 
\begin{lemma}\label{lem_unbounded_A}
Consider Assumption~\ref{ass:operators} and define the (unbounded) operator 
\begin{equation}
	\label{eqn_unbounded_A}
	A_\gamma := \begin{bmatrix}
	- \gamma \A & \K^*\\
	- \K & -\D/\gamma
	\end{bmatrix} 
	\colon D(A_\gamma)\subseteq (\calH_{\ker} \times \M) \to \calH_{\ker} \times \M
\end{equation} 
for an arbitrary positive parameter~$\gamma>0$. 
Then, $A_\gamma$ generates a $C_0$-semigroup with the domain
\[ 
D(A_\gamma) 
= \Pker \times  \big\{ m \in \M\ |\ \exists\, \mstar \in \calH_{\ker} \colon (\mstar, q_{\ker})_{\calH} 
= \langle \K^* m, q_{\ker}\rangle\ \text{for all } q_{\ker} \in \Pker \big\}.
\]
\end{lemma}
\begin{proof}
The proof is given in Appendix~\ref{appendix}. 
\end{proof}
With the previous lemmata, we are now in the position to discuss the unique solvability of the PDAE~\eqref{eq:PDAE:eps}. 
%
%
\subsection{Existence of mild and classical solutions}\label{sect:existence:mildClassical}
In this subsection, we first discuss the existence of mild solutions and turn to classical solutions afterwards. We emphasize that the property of~$\eps$ being small is not needed for the here presented existence results. 
\begin{proposition}[Existence of a mild solution]
\label{prop_mild_sol_B}
Consider Assumption~\ref{ass:operators} and right-hand sides $\g=\g_1 + \g_2$ with $\g_1\in H^1(0,T;\P^*)$, $\g_2\in L^2(0,T;\cH^*)$, $\f \in L^2(0,T;\M^*)$, and $\h \in H^1(0,T;\Q^*)$. Further assume initial data $p(0) \in \cH$ with $p(0)-\B^-\h(0)\in\cHker$ and $m(0)\in \M$. In this case, there exists a unique mild solution $(p,m,\lambda)$ of~\eqref{eq:PDAE:eps} with
\[
	p\in C([0,T],\cH) \cap H^1(0,T;\Pker^*) 
	\qquad \text{and} \qquad 
	m\in C([0,T],\M).
\]
Moreover, the Lagrange multiplier~$\lambda$ exists in a distributional sense with a regular primitive in the space~$C([0,T],\Q)$ and it satisfies that 
\[ 
	\dot{p} + \B^* \lambda \in L^2(0,T;\P^*).
\]
\end{proposition}
\begin{proof}
Let $(\pb, \mb, \lb) \in H^1(0,T;\Pker) \times H^1(0,T;\M) \times  H^1(0,T;\Q)$ be the unique solution of system~\eqref{eq:op:stationary} with given right-hand side~$\g_1$,  cf.~Lemma~\ref{lem_op_B_station}. 
The introduction of 
\begin{align}
\label{def_tildep}
	\pt := p - \pb - \B^- h,\qquad 
	\mt := m - \mb, \qquad
	\lt := \lambda - \lb
\end{align}
leads, together with~\eqref{eq:PDAE:eps}, to the system 
\begin{alignat*}{5}
	\dot \pt &\ +\ & \A\pt &\ -\ & \K^* \mt &\ +\ & \B^*\lt &= \g_2 + C_{\A} \pb - \dot \pb - \A\B^- \h - \B^- \dot \h  && \qquad \text{in }\P^*, \\
	\eps\, \dot \mt &\ +\ & \K \pt &\ +\ & \D \mt &  &  &= \f - \K \B^- \h - \eps\, \dot \mb  && \qquad \text{in }\M^*, \\
	& & \B \pt & & & & &= 0 && \qquad \text{in } \Q^*
\end{alignat*}
with initial values $\pt(0) = p(0)-\pb(0) - \B^- h(0) \in \cHker$ and $\mt(0) = m(0) - \mb(0) \in \M$. 
Since we have~$\B \pt = 0$, the solution~$\pt$ takes values in~$\Pker$. Hence, we can reduce the system to an unconstrained problem, for which we can prove existence of a solution, cf.~the details in Appendix~\ref{appendix}. 
\end{proof}
Considering higher regularity for the given data, we can show the existence of a classical solution. For this, we again analyze the corresponding Cauchy problem. 
\begin{proposition}[Existence of a classical solution]
\label{prop_clas_sol_B}
Let Assumption~\ref{ass:operators} hold and the right-hand sides satisfy $\g=\g_1 + \g_2$ with $\g_1\in H^2(0,T;\P^*)$, $\g_2\in H^1(0,T;\cH^*)$, $\f \in H^1(0,T;\M^*)$, and $\h \in H^2(0,T;\Q^*)$.  
Further assume consistent initial data $p(0) \in \P$, i.e., $\B p(0)=\h(0)$, $m(0)\in \M$, and the existence of an element $\mstar \in \cHker$ with 
\begin{align}
\label{def_mstar}
	(\mstar, q_{\ker})_{\cH} 
	= \langle \K^* m(0) + \g_1(0),\, q_{\ker}\rangle
\end{align}
for all $q_{\ker}\in \Pker$. Then there exists a unique classical solution $(p,m,\lambda)$ of~\eqref{eq:PDAE:eps} with
\begin{align*}
	p\in C([0,T],\P) \cap C^1([0,T],\cH),\qquad  
	m\in C^1([0,T],\M),\qquad
	\lambda \in C([0,T], \Q).
\end{align*}
\end{proposition}
\begin{proof}
The proof is given in Appendix~\ref{appendix}. 
\end{proof}
\begin{remark}
\label{rem:m}
Exemplarily, we comment on the existence of~$\mstar$ in the setting of Example~\ref{exp:pressure}, cf.~Section~\ref{sect:PDAE:example}. Here, a sufficient condition for the existence of~$\mstar \in \cHker$ is $m(0)\in H^1(0,1)$ and~$\g_1(0) \in L^2(0,1)$. The function~$\mstar$ then equals the spatial derivative of $m(0)$ plus~$\g_1(0)$. 
\end{remark}
%
%
%
\subsection{An energy estimate}\label{sect:existence:energy}
In the upcoming analysis of Section~\ref{sect:expansion}, we need energy estimates of the mild solution of~\eqref{eq:PDAE:eps}. 
One particular estimate is subject of the following lemma. 	
\begin{lemma}
Suppose that the assumptions of Proposition~\ref{prop_mild_sol_B} are satisfied. Then, the mild solution fulfills the energy estimate 
\begin{align}
\Vert p(t)\Vert^2_\cH + \eps\, \Vert m(t) \Vert_\M^2 + c_\D \int_0^t\Vert m (s)\Vert_\M^2 \ds 
\lesssim e^{(1+2C_{\A})t}\, \big( \Vert p(0) \Vert^2_\cH + \eps\, \Vert m(0)\Vert^2_\M + C^2_\text{data} \big) 
\label{eqn_estimate_clas}
\end{align}
with the constant 
\begin{align}
\label{def:Cdata}
C_\text{data}^2
:= \|\g_1\|^2_{H^1(0,T;\P^*)} + \|\g_2\|^2_{L^2(0,T;\cH^*)} + \|\f\|^2_{L^2(0,T;\M^*)} + \|\h\|^2_{H^1(0,T;\Q^*)}.
\end{align}
\end{lemma}
\begin{proof}
We first consider the classical solution from Proposition~\ref{prop_clas_sol_B} and $(\pt, \mt)$ as defined in~\eqref{def_tildep}.  
As shown in~Appendix~\ref{appendix}, this pair solves system~\eqref{eqn_inproof_ph} such that the application of $\pt$ and $\mt$ as test functions yields 
\begin{align*}
&\Vert \pt(t) \Vert^2_\cH + \eps\, \Vert \mt(t) \Vert_\M^2 + c_\D \int_0^t\Vert \mt (s)\Vert_\M^2 \ds \\
&\quad \le \Vert \pt(0) \Vert^2_\cH + \eps\, \Vert \mt(0)\Vert^2_\M 
+ (1+2C_{\A})\int_0^t \| \tilde p(s) \|^2_\cH  \ds \\
&\qquad + \int_0^t \big\| (\g_2 + C_{\A}\pb- \dot{\pb} - \A\B^- \h - \B^-\dot\h)(s) \big\|^2_{\cH^*}  
+ \tfrac{1}{c_\D} \big\| (\f - \K\B^-\h - \eps\, \dot{\mb})(s) \big\|^2_{\M^*} \ds. 
\end{align*}
By the continuity result of Lemma~\ref{lem_op_B_station}, the last two integrals are bounded by a multiple of~$C_\text{data}^2$. Hence, an application of Gronwalls lemma hence yields 
\begin{align*}
\Vert \pt(t) \Vert^2_\cH + \eps\, \Vert \mt(t) \Vert_\M^2 + c_\D \int_0^t\Vert \mt (s)\Vert_\M^2 \ds 
\lesssim e^{(1+2C_{\A})t}\, \big( \Vert p(0) \Vert^2_\cH + \eps\, \Vert m(0)\Vert^2_\M + C^2_\text{data} \big).
\end{align*}
Using the continuity of the semigroup generated by $A_\gamma$, the density of $D(A_\gamma)$ in $\cHker \times \M$, and the density of $H^{\ell+1}(0,T;X)$ in $H^{\ell}(0,T;X)$ for a Banach space $X$ and $\ell\ge 0$, it follows that this estimate is also satisfied for~$(\pt, \mt)$ being the mild solution of~\eqref{eqn_inproof_ph}. 
In particular, this means that the latter energy estimate is already valid under the assumptions of Proposition~\ref{prop_mild_sol_B}. 
Finally, by~$p = \pt + \pb + \B^- h$ and $m = \mt + \mb$ we conclude (again by Lemma~\ref{lem_op_B_station}) the claimed estimate~\eqref{eqn_estimate_clas}. 
\end{proof}
%
%
\subsection{Parabolic limit case}\label{sect:existence:noEps}
Finally, we consider the limit case of the PDAE~\eqref{eq:PDAE:eps} for $\eps=0$, which leads to a PDAE of parabolic type. 
The corresponding solution is denoted by~$(p_0, m_0, \lambda_0)$ and fulfills the system equations 
\begin{subequations}
	\label{eq:PDAE:noEps}
	\begin{alignat}{5}
	\dot p_0&\ +\ & \A p_0&\ -\ &\K^* m_0&\ +\ & \B^*\lambda_0 &= \g &&\qquad \text{in }\P^*, \label{eq:PDAE:noEps:a}\\
	& & \K p_0&\ +\ & \D m_0& & &= \f  &&\qquad \text{in }\M^*, \label{eq:PDAE:noEps:b}\\
	& & \B p_0& & & &  &= \h &&\qquad \text{in } \Q^*. \label{eq:PDAE:noEps:c}
	\end{alignat}
\end{subequations}
Since there is only a single differential variable (in time) left -- recall that~$\dot{m}_0$ does not appear anymore -- we consider as initial condition~$p_0(0) = p(0)$. As before, we discuss the existence of solutions for different regularity assumptions. 
\begin{proposition}[Existence of a weak solution $(p_0,m_0)$]
\label{prop_p0m0}
Consider Assumption~\ref{ass:operators} and right-hand sides $\g\in L^2(0,T;\P^*)$, $\f \in L^2(0,T;\M^*)$, and $\h \in H^1(0,T;\Q^*)$. Further assume that the initial data is consistent in the sense that~$p_0(0)-\B^-\h(0)\in\cHker$. Then, system~\eqref{eq:PDAE:noEps} has a unique weak solution with
\[
	p_0 \in L^2(0,T;\P) \cap C([0,T], \cH)
	\qquad\text{and}\qquad 
	m_0 \in L^2(0,T;\M).
\]
The Lagrange multiplier~$\lambda_0$ exists in a distributional sense with $\dot{p}_0 + \B^*\lambda_0 \in L^2(0,T;\P^*)$.
\end{proposition}
\begin{proof}
Since equation~\eqref{eq:PDAE:noEps:b} is stated in $\M^* \cong\M$ and~$\D$ is invertible, we can insert this equation into~\eqref{eq:PDAE:noEps:a}, which results in 
\begin{subequations}
	\label{eq:PDAE:noEps:red}
	\begin{alignat}{5}
		\dot p_0 + (\wlap + \A) p_0& + \B^*\lambda_0 &\ =\ & \g +  \K^* \D^{-1} \f &&\qquad \text{in }\P^*,\label{eq:PDAE:noEps:red:a}\\
		\B p_0&  &\ =\ & \h &&\qquad \text{in } \Q* \label{eq:PDAE:noEps:red:c}
	\end{alignat}
\end{subequations}
with the operator $\wlap$ introduced in Lemma~\ref{lem_wlap}. Note that the operator~$\wlap + \A$ satisfies a G\aa{}rding inequality on $\Pker$. As a result, the existence of a unique partial solution $(p_0,\lambda_0)$ with the claimed properties follows by~\cite[Th.~3.3]{EmmM13}. 
Finally, with equation~\eqref{eq:PDAE:noEps:b} the flow variable is given by~$m_0 = \D^{-1}(\f - \K p_0)$ and therefore unique and an element of $L^2(0,T;\M)$.
\end{proof}
For the subsequent analysis we also need solutions with higher regularity including a continuous multiplier~$\lambda_0$. This is subject of the following proposition. 
\begin{proposition}[Weak solution with higher regularity]
\label{prop_weak_op_B_noeps_reg} 
Consider Assumption~\ref{ass:operators} and right-hand sides $\f \in H^1(0,T;\M)$, $\g \in W(0,T; \P^*, \Pker^*)$, and $\h \in H^2(0,T;\Q^*)$. Further assume consistent initial data $p_0(0) \in \P$ with $\B p_0(0) = \h(0)$ and the existence of a function $\pstar\in \cHker$ such that 
\begin{align}
\label{def_pstar}
	(\pstar, q_{\ker})_{\calH} 
	= \big\langle \g(0) - (\wlap + \A) p_0(0) - \B^- \dot h (0) + \K^* \D^{-1}\f(0), q_{\ker} \big\rangle
\end{align}
for all $q_{\ker}\in \Pker$. 
Then the solution of system~\eqref{eq:PDAE:noEps} satisfies 
\[
	p_0 \in  H^1(0,T; \P) \cap C^1([0,T],\cH), \qquad
	m_0 \in  H^1(0,T; \M), \qquad
	\lambda_0 \in C([0,T],\Q).
\]
\end{proposition}
\begin{proof}
Due to~\eqref{eq:PDAE:noEps:red} the function $\pt_0 := p_0 -\B^-\h$ satisfies $\B\pt_0 = 0$ and  
\begin{equation}
\label{eq:PDAE:noEps_pot}
	\dot \pt_0 + (\wlap+\A) \pt_0 
	= \g - (\wlap + \A) \B^-\h - \B^-\dot \h + \K^* \D^{-1} \f \qquad \text{in }\Pker^*.
\end{equation}
The given assumptions imply that the right-hand side is an element of $H^1(0,T;\Pker^*)$ and~$\dot \pt_0(0) = \pstar \in \cHker$. 
By \cite[Th.~IV.27.2]{Wlo92} we conclude that $\pot \in H^1(0,T;\Pker)\cap C^1([0,T],\cH)$. This, in turn, implies the stated regularity for~$m_0 = \D^{-1}(\f - \K p_0)$ and the existence of~$\lambda_0$ by~\cite[Lem.~III.4.2]{Bra07}. 
\end{proof}
\begin{remark}
\label{rem:p}
Considering again Example~\ref{exp:pressure}, a sufficient condition for the existence of~$\pstar$ as in~\eqref{def_pstar} is that the difference between $\f(0)$ and the spatial derivative of the initial value $p(0)$ has a weak derivative in space and that $\g(0)$ is an $L^2$-function. Hence, for sufficiently smooth right-hand sides this reduces to~$p_0(0)\in H^2(\Omega)$.  
\end{remark}
\section{Expansion of the Solution}\label{sect:expansion}
This section is devoted to an expansion of the solution triple $(p,m,\lambda)$ of~\eqref{eq:PDAE:eps} in terms of the small parameter~$\eps$. Assuming $\eps \ll 1$, we consider the expansion 
\begin{align}
  p = p_0 + \eps\, p_1 + \dots, \qquad
  m = m_0 + \eps\, m_1 + \dots, \qquad
  \lambda = \lambda_0 + \eps\,\lambda_1 + \dots. \label{eqn_epsExp}
\end{align}
Note that the triple $(p_0, m_0, \lambda_0)$ is the solution of~\eqref{eq:PDAE:eps} in the limit case~$\eps=0$, i.e., of system~\eqref{eq:PDAE:noEps}. The aim of this section is to prove approximation properties of first order for~$p_0$ and $m_0$ as well as of second order for  
\[
  \hat p := p_0 + \eps p_1 
  \qquad\text{and}\qquad
  \hat m := m_0 + \eps m_1.
\]
Resulting approximation properties of the Lagrange multiplier will then be discussed in Section~\ref{sect:expansion:lagrange}. 
%
%
\subsection{First-order approximation}\label{sect:expansion:first}
We first discuss the approximation property of the pair $(p_0,m_0)$, i.e., we compare the solutions of the two systems~\eqref{eq:PDAE:eps} and~\eqref{eq:PDAE:noEps}. For the third example discussed in Section~\ref{sect:PDAE:example}, it was shown in~\cite[Th.~1]{EggK17ppt} that this approximation is of order $\sqrt\eps$ and -- under certain assumptions on the initial data -- of order $\eps$. We will rediscover this result in the more general setting, using an alternative technique of proof. For this, we consider the difference of the two systems~\eqref{eq:PDAE:eps} and~\eqref{eq:PDAE:noEps}, which leads to 
\begin{subequations}
\label{eqn_inproof_pp0}
\begin{alignat}{5}
 	\ddt (p-p_0)&\,+\,& \A (p-p_0)&\,-\,& \K^* (m-m_0) & + \B^*(\lambda-\lambda_0)\ &=& \ 0  &&\qquad \text{in } \P^*, \label{eqn_inproof_pp0_a}\\
	&  & \K (p-p_0)&\,+\,& \D (m-m_0) & &=&\ -\eps\, \dot m &&\qquad \text{in }\M^*, \label{eqn_inproof_pp0_b} \\	
	& & \B (p-p_0) & & & &=&\ 0 &&\qquad \text{in } \Q^*. \label{eqn_inproof_pp0_c}
\end{alignat}
\end{subequations}
The corresponding initial condition reads $(p-p_0)(0) = 0$. 
Due to~\eqref{eqn_inproof_pp0_c} we expect the difference $p-p_0$ to take values in~$\Pker$. Because of this, we will often restrict the test space in~\eqref{eqn_inproof_pp0_a} to~$\Pker$. In this case, the Lagrange multipliers vanish and we have the equation  
\begin{align}
\label{eqn_inproof_pp0_a_ker}
  \ddt (p-p_0) + \A (p-p_0) - \K^* (m-m_0) = 0  \qquad \text{in } \Pker^*.
\end{align}
In the following, we consider two approaches to derive error estimates: First, we understand~\eqref{eqn_inproof_pp0} as a parabolic system with a right-hand side $-\eps\,\dot m$. Second, we will subtract~$\eps\, \dot m_0$ from~\eqref{eqn_inproof_pp0_b} and consider the result as a hyperbolic system with a right-hand side $-\eps\,\dot m_0$.  
\begin{theorem}[First-order approximation I]
\label{thm_p0m0}
Suppose that all assumptions of Proposition~\ref{prop_clas_sol_B} are satisfied such that the PDAE~\eqref{eq:PDAE:eps} has a classical solution. Then the difference of $(p,m)$ and $(p_0, m_0)$ is bounded for $0\le t\le T$ by
\begin{multline*}
  \Vert (p-p_0)(t) \Vert_\cH^2 + \cLap \int_0^t \Vert (p-p_0)(s) \Vert^2_\P \ds  + c_\D \int_0^t \Vert (m-m_0)(s) \Vert_\M^2 \ds \\
  \lesssim \eps\, e^{(1+4C_{\A})t}\, \Vert \f(0) - \K p(0) - \D m(0) \Vert^2_{\M^*}
  + \eps^2\, e^{(1+4C_{\A})t}\, \widetilde C_\text{data}
\end{multline*}
with a constant~$\widetilde C_\text{data}$ only depending on~$\Vert p(0) \Vert_\cH$, $\|\mstar\|_{\cH}$, and the right-hand sides~$\|\g_1\|_{H^2(0,T;\P^*)}$, $\|\g_2\|_{H^1(0,T;\cH^*)}$, $\|\f\|_{H^1(0,T;\M^*)}$, and $\|\h\|_{H^2(0,T;\Q^*)}$.
\end{theorem}
\begin{proof}
We test \eqref{eqn_inproof_pp0_a_ker} and \eqref{eqn_inproof_pp0_b} by $p-p_0$ and $m-m_0$, respectively. Adding and integrating the resulting equations, we obtain by Young's and Gronwall's inequality
\begin{equation}
\label{eqn_inproof_p0m0}
  \Vert (p-p_0)(t) \Vert^2_\cH + c_\D \int_0^t \Vert (m-m_0)(s) \Vert_\M^2 \ds \leq 
  \frac{\eps^2}{c_\D} e^{2C_{\A}t}\int_0^t \Vert \dot m(s)\Vert_\M^2 \ds.
\end{equation}
On the other hand, considering test functions $(\D^{-1})^* \K(p-p_0)$ in place of~$m-m_0$, we conclude by the continuity of the operators that  
\begin{equation}
\label{eqn_inproof_p0m0_2}
  \Vert (p-p_0)(t) \Vert^2_\cH + \cLap \int_0^t \Vert (p-p_0)(s) \Vert^2_{\P} \ds
    \le 
    \eps^2 \frac{ C^2_\K}{\cLap c_\D^2} e^{2C_{\A}t} \int_0^t \Vert \dot m(s)\Vert^2_\M \ds.
\end{equation}
Note that we have used that~$\dot m(t)\in\M$, which is guaranteed by Proposition~\ref{prop_clas_sol_B}. 
Further, we note that the derivative of the classical solution~$(m,p,\lambda)$ is again a mild solution of~\eqref{eq:PDAE:eps} where we replace the right-hand sides by their temporal derivatives. Hence, we can apply estimate~\eqref{eqn_estimate_clas}, which yields 
\begin{align*}
  c_\D \int_0^t\Vert \dot m (s)\Vert^2_\M \ds
  \lesssim  e^{(1+2C_{\A})t}\, \big( \Vert \dot p(0) \Vert^2_\cH + \eps\, \Vert \dot m(0)\Vert^2_\M +  \dot C^2_\text{data} \big)
\end{align*}
with the constant~$\dot C^2_\text{data}$ defined accordingly to~\eqref{def:Cdata}, namely 
\[
  \dot C_\text{data}^2
  := \|\g_1\|^2_{H^2(0,T;\P^*)} + \|\g_2\|^2_{H^1(0,T;\cH^*)} + \|\f\|^2_{H^1(0,T;\M^*)} + \|\h\|^2_{H^2(0,T;\Q^*)}. 
\] 
It remains to bound the initial values of $\dot p$ and $\dot m$. 
For $\dot p(0)$ we use \eqref{eq:PDAE:eps:a} and the fact that~$\K^* m(0) + \g_1(0)$ has a representation in~$\cHker$ by~\eqref{def_mstar},
\[
  \|\dot p(0)\|_{\cH}^2
  = \langle \g(0) - \A p(0) + \K^* m(0), \dot p(0)\rangle    
  = \langle \mstar + \g_2(0) - \A p(0), \dot p(0)\rangle.
\]
Thus, we have~$\|\dot p(0)\|_{\cH}^2 \lesssim \|\mstar\|^2_{\cH} + \|\g_2\|^2_{H^1(0,T;\cH^*)} + \|p(0)\|_{\cH}^2$. 
For an estimate of $\dot m(0)$, we simply apply~\eqref{eq:PDAE:eps:b} to obtain
\begin{align*}
  \eps\, \Vert \dot m(0)\Vert^2_\M
  = \eps^{-1}\, \Vert \f(0) - \K p(0) - \D m(0) \Vert^2_{\M^*}.
\end{align*}
In total, this gives

\begin{equation*}
  \eps^2 \int_0^t\Vert \dot m (s)\Vert^2_\M \ds
  \lesssim e^{(1+2C_{\A})t} \Big[ \eps\, \Vert \f(0) - \K p(0) - \D m(0) \Vert^2_{\M^*}
  + \eps^2\, \big( \Vert p(0) \Vert^2_\cH + \|\mstar\|^2_{\cH} + \dot C^2_\text{data}\big) \Big],
\end{equation*}
which completes the proof.  
\end{proof}
As mentioned above, the second approach considers system~\eqref{eqn_inproof_pp0} as a hyperbolic system. For this, we need to assume the existence of $\dot m_0$ (cf.~Proposition~\ref{prop_weak_op_B_noeps_reg}), which then appears on the right-hand side. 
\begin{theorem}[First-order approximation II]
\label{thm_p0m0_2}
Under the assumptions of Proposition~\ref{prop_weak_op_B_noeps_reg}, the difference of the solutions of~\eqref{eq:PDAE:eps} and~\eqref{eq:PDAE:noEps} satisfy for~$0\le t\le T$ the estimate 
\begin{multline*}
  \Vert (p-p_0)(t) \Vert_\cH^2 + \eps\, \Vert (m-m_0)(t) \Vert_\M^2 + c_\D \int_0^t \Vert (m-m_0)(s) \Vert_\M^2 \ds\\
  \lesssim \eps\, e^{2C_{\A}t}\, \Vert m(0)-m_0(0) \Vert^2_{\M}
  + \eps^2\, e^{4C_{\A}t}\, \widetilde C_\text{data}
\end{multline*}
with a constant~$\widetilde C_\text{data}$ only depending on the initial data, $\Vert \pstar \Vert^2_{\cHker}$, and the right-hand sides $\Vert \f\Vert^2_{H^1(0,T;\M)}$, $\Vert \g\Vert^2_{W(0,T;\P^\ast,\Pker^\ast)}$, and $\Vert \h\Vert^2_{H^2(0,T;\Q^\ast)}$. 
\end{theorem}
\begin{proof}
Adding $-\eps\, \dot m_0$ to~\eqref{eqn_inproof_pp0_b}, we obtain the to~\eqref{eqn_inproof_pp0} equivalent system 
\begin{alignat*}{5}
 \ddt (p-p_0)&\, +\, & \A (p-p_0)&\,-\,& \K^* (m-m_0) & + \B^*(\lambda-\lambda_0)\ &=&\ 0  &&\qquad \text{in }\P^*, \\ 
 \eps\, \ddt (m-m_0)&\, +\, &\K (p-p_0)&\,+\,& \D (m-m_0) & &=&\ -\eps\, \dot m_0 &&\qquad \text{in }\M^*, \\ 
 & & \B (p-p_0) & & & &=&\ 0 &&\qquad \text{in } \Q^*.
\end{alignat*}
We would like to emphasize that this system has the same structure as the original PDAE~\eqref{eq:PDAE:eps}. 
By Proposition~\ref{prop_weak_op_B_noeps_reg} we conclude that~$\dot m_0 \in L^2(0,T;\M)$. Hence, we can apply the estimate for mild solutions~\eqref{eqn_estimate_clas} with right-hand sides~$\g=0$, $\f=-\eps\, \dot m_0$, and~$\h=0$. 
With $p(0)=p_0(0)$ we conclude 
\begin{align*}
  \Vert (p-p_0)(t) \Vert_\cH^2 + \eps\,& \Vert (m-m_0)(t) \Vert_\M^2 + c_\D \int_0^t \Vert (m-m_0)(s) \Vert_\M^2 \ds\\
  &\qquad\quad\lesssim e^{2C_{\A}t}\, \big(\eps\, \Vert m(0)-m_0(0)\Vert^2_\M + \eps^2\|\dot m_0\|^2_{L^2(0,T;\M^*)} \big). 
\end{align*}
The second term of the right-hand side is bounded by $\eps^2  \widetilde C_\text{data}$, since~$\D \dot m_0 = \dot \f - \K \dot p_0= \dot \f - \K \B^- \dot \h - \K \dot \pt_0$, where $\dot \pt_0$ solves the formal derivative of~\eqref{eq:PDAE:noEps_pot} with initial value~$\pstar$; see also~\cite[Th. IV.27.2]{Wlo92}. 
\end{proof}
\begin{remark}
\label{rem_firstorderineps}
In the case $\eps\, \dot m(0) = \f(0) -\K p(0) - \D m(0) = 0$, which is equivalent to $m(0) = m_0(0)$, Theorems~\ref{thm_p0m0} and~\ref{thm_p0m0_2} state that~$p_0$ is a first-order approximation of~$p$ in terms of~$\eps$, measured in~$L^\infty(0,T;\cH)$ and~$L^2(0,T;\P)$. 
Further, $m_0$ is a first-order approximation of~$m$ in $L^2(0,T;\M)$ and, under the conditions of Theorem~\ref{thm_p0m0_2}, an approximation of order~$\sfrac 12$ in $C([0,T],\M)$.
\end{remark}
\begin{remark}
\label{rem_finiteinfinite}
In the finite-dimensional case, the obtained results for $m-m_0$ also match with~\cite[Ch.~2.5, Th.~5.1]{KokKO99}. For~$p-p_0$, however, one has $\|p-p_0\|_{L^\infty(0,T)} = \O(\eps)$ in the finite-dimensional setting, independent of the initial data. Ghisi and Gobbino showed the same convergence rates for the damped wave equation (Example~\ref{exp:wave}) with homogeneous right-hand side and regular enough initial values as well as proved their sharpness \cite[Th.~2.2 \&~2.3]{GhiG06}. Our assumed initial values, however, are less regular and a careful analysis shows that the here derived estimates are optimal \cite[App.~C]{AltZ18ppt}. This is also numerically confirmed in Section~\ref{sect:numerics}.
\end{remark}
\begin{remark}
If the operator $\A$ is non-negative in addition to Assumption~\ref{ass:operators}, then the constant $C_{\A}$ in the estimates of Theorems~\ref{thm_p0m0} and~\ref{thm_p0m0_2} can be set to zero. If~$\A$ is elliptic, then the bounds can be further improved such that the exponential function therein is strictly monotonic decreasing. 
\end{remark}
%
%
\subsection{Second-order approximation}\label{sect:expansion:second} 
So far, we have shown that $(p_0,m_0)$ provides a first-order approximation of~$(p,m)$ if the initial value of $m$ is chosen appropriately. We now include the second term of the expansion~\eqref{eqn_epsExp}. To be precise, we analyze the approximation properties of~$\hat p = p_0 + \eps\, p_1$ and~$\hat m = m_0 + \eps\, m_1$. We are especially interested in the needed regularity assumptions to gain an additional half or full order in terms of $\eps$. 

As a first step, we note that the tuple $(p_1, m_1, \lambda_1)$ introduced in~\eqref{eqn_epsExp} solves the PDAE   
\begin{subequations}
\label{eqn_op_B_approx_1}
\begin{alignat}{5}
  \dot p_1 &\,+\, & \A p_1 &\, - \, & \K^* m_1 &\, +\, &\B^*\lambda_1 &= 0  &&\qquad \text{in }\P^*, \label{eqn_op_B_approx_1_a}\\
  & & \K p_1 &\, +\, & \D m_1 & & &= -\dot m_0  &&\qquad \text{in }\M^*, \label{eqn_op_B_approx_1_b}\\
  & & \B p_1 &  & & & &= 0 &&\qquad \text{in } \Q^* \label{eqn_op_B_approx_1_c}
\end{alignat}
\end{subequations}
with the initial condition $p_1(0) = 0$. 
For the solvability of system~\eqref{eqn_op_B_approx_1} we note that it has the same structure as~\eqref{eq:PDAE:noEps}. Hence, we only need to analyze the regularity of the right-hand side, i.e., of $\dot m_0$. 
The weak differentiability of $m_0$ has been discussed in Proposition~\ref{prop_weak_op_B_noeps_reg} such that an application of Proposition~\ref{prop_p0m0} leads to the following result.
\begin{proposition}[Existence of a weak solution $(p_1, m_1)$]
\label{prop_p1m1}
Given the assumptions of Proposition~\ref{prop_weak_op_B_noeps_reg}, system~\eqref{eqn_op_B_approx_1} is uniquely solvable with 
\[
  p_1 \in L^2(0,T;\Pker) \cap C([0,T], \cH)
  \qquad\text{and}\qquad 
  m_1 \in L^2(0,T;\M).
\]
Furthermore, $\lambda_1$ exists in a distributional sense.
\end{proposition}
In the previous subsection, we have observed that the initial data may cause a reduction in the approximation order, cf.~Remark~\ref{rem_firstorderineps}. To focus on the improvements resulting from the incorporation of $p_1$ and $m_1$, we assume in the following that $\f(0) - \K p(0) - \D m(0) = 0$, i.e., $m(0) = m_0(0)$. 
\begin{theorem}[Second-order approximation]
\label{thm_p1m1}
Consider the assumptions of Proposition~\ref{prop_clas_sol_B} with additional regularity of the form $\g=\g_1 + \g_2$ with $\g_1\in H^3(0,T;\P^*)$, $\g_2\in H^2(0,T;\cH^*)$, $\f \in H^2(0,T;\M^*)$, and $\h \in H^3(0,T;\Q^*)$. Further assume the existence of an element~$\pstar\in \Pker$  satisfying \eqref{def_pstar}, $\f(0) - \K p(0) - \D m(0) = 0$, and that $\dot \g_1(0)=0$. In this setting, the difference of the solution of~\eqref{eq:PDAE:eps} and $(\hat p, \hat m)$ is bounded by 
\begin{align*}
  \Vert (p-\hat p)(t)\Vert^2_\cH + \cLap \int_0^t &\Vert (p-\hat p)(s) \Vert^2_{\P} \ds + c_\D \int_0^t \Vert (m-\hat m)(s)\Vert^2_\M \ds \\
  &\qquad\lesssim \eps^3\, e^{(1+4C_{\A})t}\, \Vert \dot \f(0) - \K \dot p (0)\Vert^2_{\M^*} + \eps^4\, e^{(1+4C_{\A})t}\,  \hat C_\text{data} 
\end{align*}
for $0 \le t\le T$ and with~$\hat C_\text{data}$ depending on the initial data, $\Vert \pstar\Vert^2_{\cHker}$, and the right-hand sides $\Vert \g_1\Vert^2_{H^3(0,T;\P^\ast)}$, $\Vert \g_2\Vert^2_{H^2(0,T;\cH^\ast)}$, $\Vert \f\Vert^2_{H^2(0,T;\M^\ast)}$, and~$\Vert \h\Vert^2_{H^3(0,T;\Q^\ast)}$. 
\end{theorem}
\begin{proof}
The difference of the solution~$(p, m, \lambda)$ of~\eqref{eq:PDAE:eps} and $(\hat p, \hat m, \hat \lambda)$ solves 
\begin{alignat*}{5}
\ddt (p-\hat p)&\, +\, & \A (p-\hat p)&\,-\,& \K^* (m-\hat m) & + \B^*(\lambda-\hat\lambda)\ &=&\ 0  &&\qquad \text{in }\P^*, \\ 
&  &\K (p-\hat p)&\,+\,& \D (m-\hat m) & &=&\ -\eps\, (\dot m-\dot m_0) &&\qquad \text{in }\M^*, \\ 
& & \B (p-\hat p) & & & &=&\ 0 &&\qquad \text{in } \Q^* 
\end{alignat*}
with initial condition $(p-\hat p)(0) = 0$. 
Considering $p-\hat p$ as test function in the first equation, the Lagrange multipliers vanish. Thus, following the proof of Theorem~\ref{thm_p0m0}, we obtain
\begin{multline*} 
  \Vert (p-\hat p)(t)\Vert^2_\cH + \int_0^t \Vert (p-\hat p)(s) \Vert^2_{\P} 
  + \Vert (m-\hat m)(s)\Vert^2_\M \ds
  \lesssim \eps^2\, e^{2C_{\A}t} \int_0^t \Vert (\dot m-\dot m_0)(s)\Vert^2_\M \ds
\end{multline*}
and it remains to find an estimate of the integral of the error~$\dot m-\dot m_0$. 
For this, we consider the formal derivative of system~\eqref{eqn_inproof_pp0}. Similar to the estimate~\eqref{eqn_inproof_p0m0}, we can show that 
\[
  \int_0^t \Vert (\dot m - \dot m_0)(s) \Vert_\M^2 \ds 
  \lesssim \Vert (\dot p - \dot p_0)(0)\Vert^2_\cH + \eps^2\, \int_0^t \Vert \ddot m(s)\Vert_\M^2 \ds.
\]
As mentioned in Remark~\ref{rem_firstorderineps}, the assumption $\f(0) -\K p(0) - \D m(0) = 0$ implies that $m(0)=m_0(0)$. Equation~\eqref{eqn_inproof_pp0_a} thus implies~$\dot p(0) = \dot p_0(0)$ in~$\cH$.
For an estimate of~$\ddot m$, we consider the formal derivative of~\eqref{eq:PDAE:eps} and obtain by estimate~\eqref{eqn_estimate_clas} that 
\[
  \int_0^t\Vert \ddot m (s)\Vert_\M^2 \ds 
  \lesssim e^{(1+2C_{\A})t}\, \big( \Vert \ddot p(0) \Vert^2_\cH + \eps\, \Vert\ddot m(0)\Vert^2_\M + \ddot C^2_\text{data} \big). 
\]
Here, $\ddot C^2_\text{data}$ denotes the constant similarly defined as in~\eqref{def:Cdata} but with two additional derivatives on each occurring function. For an estimate of $\ddot m(0)$ we note that~\eqref{eq:PDAE:eps:b} implies together with $\dot m(0)=0$ that 
\begin{align*}
  \eps^5\, \Vert \ddot m(0)\Vert^2_\M
  = \eps^{3}\, \Vert \dot \f(0) - \K \dot p(0) - \D \dot m(0) \Vert^2_{\M^*}
  = \eps^{3}\, \Vert \dot \f(0) - \K \dot p(0) \Vert^2_{\M^*}. 
\end{align*}
Note that the right-hand side is bounded due to~$\f \in H^2(0,T;\M^\ast)$ and $\dot p(0) = \pb(0) + \pstar + \B^-\h(0) \in \P$.

In order to find a bound of~$\ddot{p}(0)$, we apply the decomposition~$p = \pb + \pt + \B^-\h$ with~$\pb$ being the solution of~\eqref{eq:op:stationary} with right-hand side~$\g_1$ and~$\pt$ the function introduced in~\eqref{def_tildep}. 
By Lemma~\ref{lem_op_B_station} we know that~$\ddot \pb(0)$ is bounded by the~$H^3(0,T;\P^*)$-norm of $\g_1$. Further, $\B^-\ddot\h$ is bounded by the~$H^3(0,T;\Q^*)$-norm of $\h$. 
Finally, \eqref{eqn_inproof_ph_a} implies that 
\[
  \big\Vert \ddot\pt(0) \big\Vert^2_\cH
  = \big\Vert \dot\g_2(0) - \A \dot \pt(0) - \A\B^- \dot\h(0) - \ddot \pb(0) - \B^-\ddot\h(0) \big\Vert^2_{\cH^*}.
\]
Here, we used that~$\dot\mt(0)=\dot m(0) + \dot \mb(0)$ vanishes by the assumptions $\eps\, \dot m(0) = \f(0) - \K p(0) - \D m(0) = 0$ and $\dot \g_1(0)=0$.   
Finally, since $p(0)=p_0(0)$, $m(0)=m_0(0)$, and $\dot m(0)=0$ imply that~$\dot  p(0) = \dot p_0(0)$ by~\eqref{eqn_inproof_ph} and~\eqref{eq:PDAE:noEps_pot}, we  have $\dot \pt(0) = \pstar - \dot \pb(0) = \pstar \in \cHker$ and thus the term $\A \dot \pt(0)$ is bounded in $\calH^\ast$.
\end{proof}
\begin{remark}
In Theorem~\ref{thm_p1m1} we have assumed $\dot \g_1(0)=0$ to show that $\dot\mt(0)$ vanishes and, hence, $\K \dot\mt(0)$ is bounded in $\cHker^\ast$. To conclude this boundedness, however, it is sufficient (and necessary) that $\dot \g_1(0)\in \cHker^\ast$. 
\end{remark}
\begin{remark}
Under the additional assumption $0=\dot{f}(0) - \K \dot{p}(0) =\eps\, \ddot m(0)$, which is equivalent to $m(0)=\hat m(0)$ and $\dot{m}(0)=\dot{m}_0(0)$, Theorem~\ref{thm_p1m1} states that~$\hat p$ is a second-order approximation of $p$ in terms of $\eps$, measured in~$L^\infty(0,T;\cH)$ and~$L^2(0,T;\P)$, respectively. Further, $\hat m$ is a second-order approximation of $m$ measured in $L^2(0,T;\M)$. 
\end{remark}
\begin{remark}
\label{rem_finiteinfinite2}
A comparison of Theorem~\ref{thm_p1m1} with the corresponding finite-dimensional case shows once more that the results for $m-\hat m$ coincide but that the difference $p-\hat p$ converges with a half $\eps$-order less for general initial data, cf.~\cite[Ch.~2.5, Th.~5.2]{KokKO99}. For a numerical validation of this result, we refer to~Section~\ref{sect:numerics}.
\end{remark}
%
%
\subsection{Estimate of the Lagrange multiplier}\label{sect:expansion:lagrange}  
At this point, we would like to take a closer look at the Lagrange multiplier~$\lambda$ and its approximation~$\lambda_0$. Since an estimate of~$\lambda - \lambda_0$ depends on the derivatives~$\dot p$ and~$\dot p_0$, the regularity assumptions of the first-order estimates in Section~\ref{sect:expansion:first} are not sufficient. 
Hence, we consider the assumptions of Theorem~\ref{thm_p1m1}. 

By Proposition~\ref{prop_clas_sol_B} we know that~$\lambda \in C([0,T], \Q)$, whereas~$\lambda_0 \in C([0,T],\Q)$ was shown in Proposition~\ref{prop_weak_op_B_noeps_reg}.
By the inf-sup stability of~$\B$ we have 
\[
  \| \lambda(t) - \lambda_0(t) \|_\Q
  \le \frac{1}{\beta}\, \| \B^* (\lambda(t) - \lambda_0(t)) \|_{\P^*}
  = \frac{1}{\beta}\, \sup_{q\in \P} \frac{\langle \B^* (\lambda(t) - \lambda_0(t)) , q\rangle}{\| q\|_\P}
\]
and thus, by equation~\eqref{eqn_inproof_pp0_a}, 
\[
  \| \lambda(t) - \lambda_0(t) \|_\Q
  \lesssim \| \dot p(t) - \dot p_0(t)\|_{\P^*} + \| p(t) - p_0(t)\|_\cH + \| m(t)-m_0(t)\|_\M. 
\]
Hence, the $L^2(0,T;\Q)$-error of $\lambda - \lambda_0$ can be bounded by 
\begin{align*}
  &\int_0^t \| \lambda(s) - \lambda_0(s) \|_\Q^2 \ds\\
  &\quad\lesssim \int_0^t \| \dot p(s) - \dot p_0(s)\|^2_{\P^*} + \| p(s) - p_0(s)\|^2_\cH + \| m(s)-m_0(s)\|^2_\M \ds \\
  &\quad\lesssim \int_0^t \| \dot p(s) - \dot p_0(s)\|^2_{\cH} \ds + \eps^2\, e^{4C_{\A}t}\, \widetilde C_\text{data} \\
  &\quad\leq \int_0^t \frac12\, \eps^{-1/2}\, \| \dot p(s) - \dot p_0(s)\|^2_{\Pker^*} + \frac12\, \eps^{1/2}\, \| \dot p(s) - \dot p_0(s)\|^2_{\Pker} + \eps^2\, e^{4C_{\A}t}\, \widetilde C_\text{data} \\
  &\quad\lesssim \eps^{3/2}\, e^{(1+4C_{\A})t}\, \Vert \dot \f(0) - \K \dot p (0)\Vert^2_{\M^*} + (\eps^{3/2} + \eps^2)\, e^{4C_{\A}t}\, \widetilde C_\text{data} + \eps^{5/2}\, e^{(1+4C_{\A})t}\,  \hat C_\text{data}
\end{align*}
with the constants $\widetilde C_{\text{data}}$ and $\hat C_\text{data}$ from Theorems~\ref{thm_p0m0} and~\ref{thm_p1m1}, respectively. Here, we used on the one hand that $\dot p - \dot p_0$ can be estimated in $\Pker^\ast$ analogously to the lines of Theorem~\ref{thm_p0m0}. On the other hand, $\dot p - \dot p_0$ and $\dot m-\dot m_0$ satisfy the formal derivative of the parabolic system~\eqref{eqn_inproof_pp0} with the  initial condition~$\dot p(0) - \dot p_0(0)=0$. Similarly to estimate~\eqref{eqn_inproof_p0m0_2} we then have $\int_0^t \|\dot p - \dot p_0\|_{\P}^2 \ds \lesssim \eps^2 \int_0^t \|\ddot m\|^2_{\M}\ds$, where the right-hand side can be bounded in the same manner as in the proof of Theorem~\ref{thm_p1m1}.
\begin{remark}\label{rem_lambda_interpol}
If the interpolation space $[\Pker^{\hphantom{\ast}},\Pker^\ast]_{\theta}$ with $\theta>\frac{1}{2}$ can be embedded in~$\P^\ast$, then the estimate of $\int_0^t \| \lambda(s) - \lambda_0(s) \|_\Q^2 \ds$ can be improved to the order~$\O(\eps^{1+\theta})$. For more details on this, we refer the reader to~\cite[Rem.~8.39]{Zim21}.
\end{remark}
%
%
\subsection{Discussion and possible extensions}\label{sect:expansion:discussion}  
Before considering several numerical experiments, we would like to summarize the obtained results and give an outlook in view of higher-order estimates. For this, we consider sufficiently smooth data. 

In Section~\ref{sect:expansion:first}, we have shown that the general approximation of $(p,m)$ by $(p_0,m_0)$ is only of order~$\sqrt{\eps}$. However, this can be improved if the condition~$m(0) = m_0(0)$ is satisfied. This means that the reduction of the convergence rate is due to a boundary layer or inconsistent initial data. 

We would like to emphasize that the worst-case estimate does not improve if we add another term in the expansion, i.e., also the approximation by~$(\hat p, \hat m)$ is only of order~$\sqrt{\eps}$ if~$m(0) \neq m_0(0)$. 
Assuming~$m(0) = m_0(0)$, we obtain the order~$\eps^{\sfrac 32}$. This can be further improved to the full order of~$\eps^2$ if the conditions~$m(0)=\hat m(0)$ and~$\dot{m}(0)=\dot{m}_0(0)$ are satisfied. Again, we observe that the reduction of the convergence rate is caused by a boundary layer.

If one is interested in higher-order expansions in $\eps$, then the expansion terms $\eps^\ell p_\ell$ and~$\eps^\ell m_\ell$ exist under the condition that~$\ddt[\ell] p_0$ exists. Even for moderate $\ell$, however, this requires very restrictive compatibility conditions on the right-hand sides $\g$, $\f$, and $\h$ as well as on the initial value $p_0(0)=p(0)$, which are only hardly practical; cf.~\cite{Tem82} and~\eqref{eq:PDAE:noEps_pot}. In addition, one would need $m(0)=m_0(0)$ and $m^{(k)}(0)=0$, $k=1,\ldots,\ell$, in order to obtain an approximation of order $\O(\eps^{\ell+\sfrac 1 2})$. 			
\section{Numerical Validation}\label{sect:numerics}
This final section is devoted to the numerical confirmation of the convergence results of Section~\ref{sect:expansion}, including the differences in the approximation orders in the finite- and infinite-dimensional setting. For this, we consider the propagation of gas in a single pipe of unit length, cf.~Section~\ref{sect:PDAE:example} Example~3. The associated PDE in its strong form is given by
\begin{subequations}
\label{eqn_numerical_example}
\begin{alignat}{5}
 \dot p& & &\,+\,& \partial_x m &=&\ 0  &&\qquad \text{in } (0,1), \\ 
 \eps\, \dot m&\, +\, &\partial_x p &\,+\,& m &=&\ 0 && \qquad \text{in } (0,1).
\end{alignat}
\end{subequations}
Moreover, the pressure~$p$ satisfies homogeneous Dirichlet boundary conditions. In the corresponding weak formulation, these boundary conditions are included explicitly in form of a constraint. With the trace operator denoted by~$\B$, $\A:=0$, and~$\D:=\id$, this then leads to the PDAE~\eqref{eq:PDAE:eps} with vanishing right-hand sides. In this particular case, the associated Lagrange multiplier~$\lambda$ equals the trace of the mass flux~$m$ (if $m$ is sufficiently smooth).

The numerical experiments of this section illustrate the transition of the approximation order in terms of~$\eps$ from the finite- to the infinite-dimensional setting, cf.~Remarks~\ref{rem_finiteinfinite} and~\ref{rem_finiteinfinite2}. Recall that Theorems~\ref{thm_p0m0} and~\ref{thm_p0m0_2} imply the estimate 
\begin{equation*}
\|p-p_0\|_{C(0,T;\,L^2(0,1))} + \|m-m_0\|_{L^2(0,T;\,L^2(0,1))}
\le C_\text{cons} \sqrt{\eps} + \O(\eps)
\end{equation*}
with $C_\text{cons} = 0$ if the initial data satisfies~$\partial_x p(0) = -m(0)$. 
In finite dimensions, one shows~$\|p-p_0\|_{C(0,T;\,\R^n)} = \O(\eps)$ independent of the initial data; see~\cite[Ch.~2.5, Th.~5.1]{KokKO99}. We would like to emphasize that this is an asymptotic result, whereas the bounds in Theorem~\ref{thm_p0m0} and~\ref{thm_p0m0_2} are valid for all~$\eps>0$. Similar statements hold true for Theorem~\ref{thm_p1m1} and its finite-dimensional counterpart.  

For the numerical validation of the approximation orders obtained in Section~\ref{sect:expansion}, we consider system~\eqref{eqn_numerical_example} for different initial values. These are given by 
\begin{align}
\label{eqn:initialData:a}
p(x,0)= \sum_{k=1}^\infty \frac{\sin(\pi k x)}{ k^{1.55}} \in H^1_0(0,1)
\quad\text{and}\quad 
m(x,0) = 0 \in H^1(0,1)
\end{align}
considering Theorem~\ref{thm_p0m0} and 
\begin{align}
\label{eqn:initialData:b}
p(x,0)= 0 \in H^2_0(0,1)
\quad\text{and}\quad 
 m(x,0)= \pi \sum_{k=1}^\infty \frac{\cos(\pi k x)}{k^{0.55}} \in L^2(0,1)
\end{align}
for Theorem~\ref{thm_p0m0_2}. Note that both pairs do not satisfy the condition~$\partial_x p(0) = -m(0)$, i.e., $C_\text{cons} \not= 0$. 

For the spatial discretization, we consider enriched spectral finite elements given by $\P_n \coloneqq \operatorname{span} \{e^{-x},e^x,\sin(\pi x), \ldots, \sin(n\pi x)\} \subset \P = H^1(0,1)$ for the pressure variable and $\M_n \coloneqq \operatorname{span} \{1,\cos(\pi x), \ldots, \cos(n\pi x)\} \subset \M = L^2(0,1)$ for the mass flux, $n\in \N$. Since the space $\Q$ for the Lagrange multiplier is simply $\R^2$, there is no need for an additional discretization. 
To identify the approximation rate as a function of the discretization parameter $n$, we calculate for fixed~$n$ the differences $p_n(\, \cdot\, ; \eps) - p_{0,n}$, $m_n(\, \cdot\, ; \eps) - m_{0,n}$, and $\lambda_n(\, \cdot\, ; \eps) - \lambda_{0,n}$ for $\eps= 1/(8\sqrt{2^j})$, $j=1,\ldots,30$. For each error measure~$\operatorname{err}(n,\eps)$ -- we consider different variables in different norms -- we make the ansatz $\operatorname{err}(n,\eps) =C(n)\,\eps^{\alpha(n)}$. With this, we obtain an approximation of the rate by the median of the slopes of the logarithmic errors between two successive values of $\eps$. 

For the initial data mentioned above, Figures~\ref{fig_p0m0} and~\ref{fig_p0m0_2} illustrate the estimated approximation orders as a function of the discretization parameter~$n$. 
\begin{figure}
%
%
\definecolor{mycolor1}{rgb}{0.00000,0.44700,0.74100}%
\definecolor{mycolor2}{rgb}{0.85000,0.32500,0.09800}%
\definecolor{mycolor3}{rgb}{0.92900,0.69400,0.12500}%
\definecolor{mycolor4}{rgb}{0.49400,0.18400,0.55600}%
\definecolor{mycolor5}{rgb}{0.46600,0.67400,0.18800}%
\definecolor{mycolor6}{rgb}{0.30100,0.74500,0.93300}%
\definecolor{mycolor7}{rgb}{0.63500,0.07800,0.18400}%
\begin{tikzpicture}

\begin{axis}[%
width=4.4in,
height=1.9in,
at={(1.335in,0.797in)},
scale only axis,
xmode=log,
xmin=1,
xmax=20000,
xminorticks=true,
xlabel style={font=\color{white!15!black}},
xlabel={discretization parameter $n$},
ymin=0,
ymax=1.25,
ylabel style={font=\color{white!15!black}},
ylabel={approximation order $\alpha$},
axis background/.style={fill=white},
title style={font=\bfseries},
legend style={legend cell align=left, align=left, draw=white!15!black}
]
\addplot [color=mycolor1, line width = 1.2]
  table[row sep=crcr]{%
1	0.963961513643703\\
2	0.925961418014377\\
4	0.852963042815319\\
8	0.737788232984866\\
16	0.603327962450793\\
32	0.549137974145618\\
64	0.531239821036733\\
128	0.526366154914392\\
256	0.524742725714485\\
512	0.524298750061478\\
1024	0.524186656992558\\
2048	0.524158681005177\\
4096	0.524152208426896\\
8192	0.524150743281099\\
16384	0.524150396488785\\
};
\addlegendentry{$p-p_0 \text{ in } L^\infty(L^2)$}

\addplot [color=mycolor4, mark=triangle, line width = 1.2]
  table[row sep=crcr]{%
1	1.00057931187347\\
2	1.00019186432332\\
4	0.996035762677886\\
8	0.969285290871809\\
16	0.907860392484295\\
32	0.826029391853376\\
64	0.752089728876402\\
128	0.673881948605536\\
256	0.63761699335864\\
512	0.612933913865727\\
1024	0.592062040900314\\
2048	0.578713893435405\\
4096	0.569470021442397\\
8192	0.56269478522867\\
16384	0.557418213997209\\
};
\addlegendentry{$p-p_0 \text{ in } L^2(H^1)$}

%

\addplot [color=mycolor5, mark=square, line width = 1.2, dashed, mark options={solid}]
  table[row sep=crcr]{%
1	0.500718476933374\\
2	0.500064180899079\\
4	0.497418754650228\\
8	0.492861061721171\\
16	0.485788182681239\\
32	0.470359877520388\\
64	0.458725929318323\\
128	0.454569318660412\\
256	0.452380775350431\\
512	0.451748124390992\\
1024	0.451438381739647\\
2048	0.448918060212436\\
4096	0.448761405925558\\
8192	0.446226339566417\\
16384	0.446128516175107\\
};
\addlegendentry{$m-m_0 \text{ in } L^2(L^2)$}


\addplot [color=mycolor7, mark=o, line width = 1.2, dotted, mark options={solid}]
  table[row sep=crcr]{%
1	0.500718476933375\\
2	0.500064180899078\\
4	0.482275233162966\\
8	0.438321943134871\\
16	0.359619847710145\\
32	0.290854931692251\\
64	0.259669998659681\\
128	0.25229607798232\\
256	0.245553058029134\\
512	0.243797003248017\\
1024	0.230070177182878\\
2048	0.214175710282862\\
4096	0.18707720180236\\
8192	0.149183768354576\\
16384	0.116794753776793\\
};
\addlegendentry{$\lambda-\lambda_0 \text{ in } L^2(\R^2)$}

\end{axis}
\end{tikzpicture}%
	\caption{Estimate of the approximation order~$\alpha$ corresponding to~Theorem~\ref{thm_p0m0} with the initial values from~\eqref{eqn:initialData:a}.}
	\label{fig_p0m0}
\end{figure}
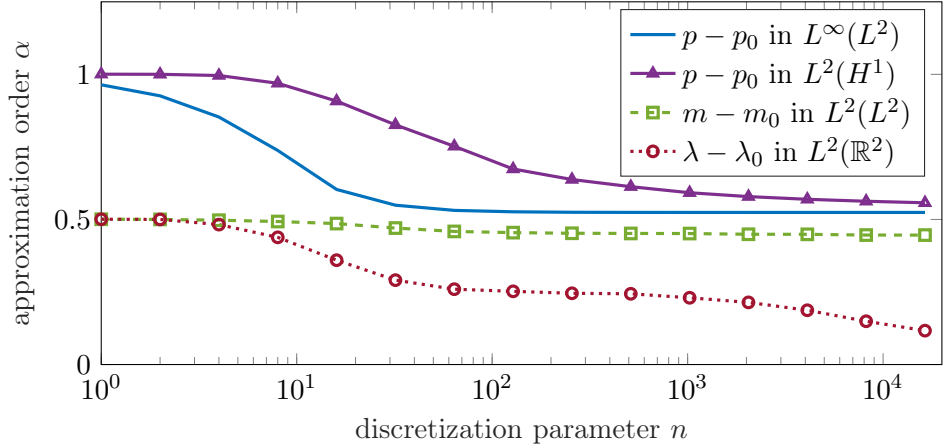%
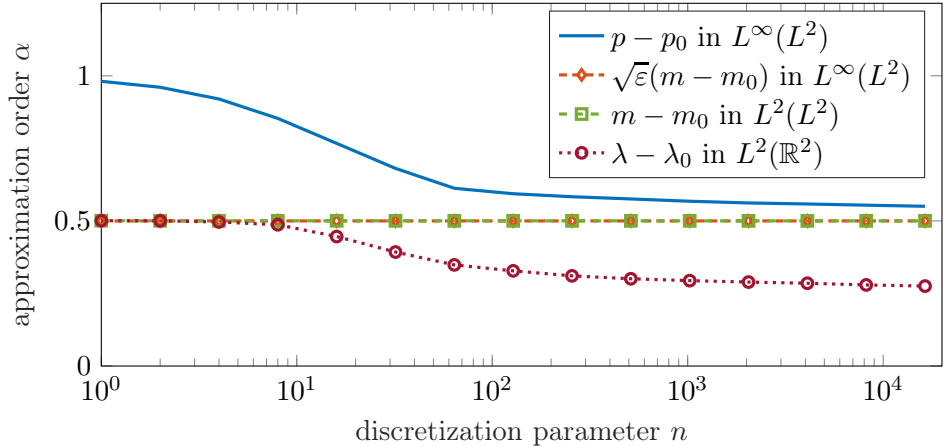
\begin{figure}
%
%
\definecolor{mycolor1}{rgb}{0.00000,0.44700,0.74100}%
\definecolor{mycolor2}{rgb}{0.85000,0.32500,0.09800}%
\definecolor{mycolor3}{rgb}{0.92900,0.69400,0.12500}%
\definecolor{mycolor4}{rgb}{0.49400,0.18400,0.55600}%
\definecolor{mycolor5}{rgb}{0.46600,0.67400,0.18800}%
\definecolor{mycolor6}{rgb}{0.30100,0.74500,0.93300}%
\definecolor{mycolor7}{rgb}{0.63500,0.07800,0.18400}%
\begin{tikzpicture}

\begin{axis}[%
width=4.4in,
height=1.9in,
at={(2.502in,1.239in)},
scale only axis,
xmode=log,
xmin=1,
xmax=20000,
xminorticks=true,
xlabel style={font=\color{white!15!black}},
xlabel={discretization parameter $n$},
ymin=0,
ymax=1.25,
ylabel style={font=\color{white!15!black}},
ylabel={approximation order $\alpha$},
axis background/.style={fill=white},
title style={font=\bfseries},
legend style={legend cell align=left, align=left, draw=white!15!black}
]
\addplot [color=mycolor1, line width = 1.2]
  table[row sep=crcr]{%
1	0.981125842248341\\
2	0.960608505423282\\
4	0.920110287806518\\
8	0.853365011140317\\
16	0.767222950795655\\
32	0.681107932535514\\
64	0.612487438206492\\
128	0.593700150747247\\
256	0.583374201859886\\
512	0.575800870483\\
1024	0.567843567266189\\
2048	0.561988491930796\\
4096	0.558368187746908\\
8192	0.554272343130492\\
16384	0.55045695959883\\
};
\addlegendentry{$p-p_0 \text{ in } L^\infty(L^2)$}

\addplot [color=mycolor2, mark=diamond, line width = 1.2, densely dashed, mark options={solid}]
  table[row sep=crcr]{%
1	0.5\\
2	0.5\\
4	0.5\\
8	0.5\\
16	0.5\\
32	0.5\\
64	0.5\\
128	0.5\\
256	0.5\\
512	0.5\\
1024	0.5\\
2048	0.5\\
4096	0.5\\
8192	0.5\\
16384	0.5\\
};
\addlegendentry{$\sqrt{\eps}(m-m_0) \text{ in } L^\infty(L^2)$}

%

\addplot [color=mycolor5, mark=square, line width = 1.2, dashed, mark options={solid}]
  table[row sep=crcr]{%
1	0.500803458057494\\
2	0.500427723740266\\
4	0.500455386511145\\
8	0.500358539950148\\
16	0.500426966637982\\
32	0.500366680350863\\
64	0.50032345685303\\
128	0.500291176097491\\
256	0.500266264363367\\
512	0.500246524136413\\
1024	0.500230541257193\\
2048	0.500217368793365\\
4096	0.500206349649652\\
8192	0.500197009365843\\
16384	0.500188932376337\\
};
\addlegendentry{$m-m_0 \text{ in } L^2(L^2)$}


\addplot [color=mycolor7, mark=o, line width = 1.2, dotted, mark options={solid}]
  table[row sep=crcr]{%
1	0.500803458057494\\
2	0.500427723740267\\
4	0.496387657756377\\
8	0.486775850747945\\
16	0.44619805843398\\
32	0.392826939835768\\
64	0.348556222175166\\
128	0.327964638937486\\
256	0.310755583390183\\
512	0.300658176495903\\
1024	0.294252549142284\\
2048	0.289328871894589\\
4096	0.285254180753105\\
8192	0.279486281263172\\
16384	0.275462753176298\\
};
\addlegendentry{$\lambda-\lambda_0 \text{ in } L^2(\R^2)$}
\end{axis}

\end{tikzpicture}%
	\caption{Estimate of the approximation order~$\alpha$ corresponding to~Theorem~\ref{thm_p0m0_2} with the initial values from~\eqref{eqn:initialData:b}.}
	\label{fig_p0m0_2}
\end{figure}%
For the pressure~$p$, we observe a fading from order one (as expected for the finite-dimensional case) to a value close to $0.5$. A rigorous analysis shows that the asymptotic limit is $0.525$. The difference in the mass flux~$m$ has a constant rate of $0.5$ for the second example (with initial data~\eqref{eqn:initialData:b}), whereas the rate is slightly lower in the first example (with initial data~\eqref{eqn:initialData:a}). We believe the latter to be a consequence of round-off errors, since analytic calculations predict a constant order of $0.5$ as well. In addition to the pressure and the mass flux, the two figures also include the estimated rate for the Lagrange multiplier~$\lambda$. Recall that we have no theoretical predictions for these particular regularity assumptions. For small $n$, it is close to $0.5$, since 
\[
\lambda_n-\lambda_{0,n} = (B_n^{\vphantom{T}}M_n^{-1}B_n^T)^{-1}B_n^{\vphantom{T}}M_n^{-1}K_n^T(m_n-m_{0,n})
\]
with the mass matrix~$M_n$, the discrete partial derivative~$K_n$, and the discrete trace operator~$B_n$. For increasing $n$, however, the approximation orders decrease. Due to the structure of system~\eqref{eqn_numerical_example}, one can prove
that $\|\partial_x (m-m_0)\|_{L^2(0,T;L^2(0,1))} = \O(1)$ for the example considered in Figure~\ref{fig_p0m0}. Since 
\[
[L^2(0,1),H^{-1}(0,1)]_{1/2-\delta} 
= [H^{1/2-\delta}(0,1)]^\ast 
\hookrightarrow [H^1(0,1)]^\ast
\]
for every $\delta\in (0,0.25)$, cf.~\cite[Ch.~1, Th.~11.1 \& 12.3]{LioM72}, the limit rate in the first example is $0.25$; see also~Section~\ref{sect:expansion:lagrange}. 
For the second example, the rate~$\alpha$ analytically tends to zero. 

\begin{figure}
%
%
\definecolor{mycolor1}{rgb}{0.00000,0.44700,0.74100}%
\definecolor{mycolor2}{rgb}{0.85000,0.32500,0.09800}%
\definecolor{mycolor3}{rgb}{0.92900,0.69400,0.12500}%
\definecolor{mycolor4}{rgb}{0.49400,0.18400,0.55600}%
\definecolor{mycolor5}{rgb}{0.46600,0.67400,0.18800}%
\definecolor{mycolor6}{rgb}{0.30100,0.74500,0.93300}%
\definecolor{mycolor7}{rgb}{0.63500,0.07800,0.18400}%
\begin{tikzpicture}

\begin{axis}[%
width=4.4in,
height=1.9in,
at={(2.502in,1.239in)},
scale only axis,
xmode=log,
xmin=1,
xmax=20000,
xminorticks=true,
xlabel style={font=\color{white!15!black}},
xlabel={discretization parameter $n$},
ymin=0.75,
ymax=1.55,
ylabel style={font=\color{white!15!black}},
ylabel={approximation order $\alpha$},
axis background/.style={fill=white},
title style={font=\bfseries},
legend style={legend cell align=left, align=left, draw=white!15!black}
]
\addplot [color=mycolor1, line width = 1.2]
  table[row sep=crcr]{%
1	1.00008498395309\\
2	1.00043123980669\\
4	1.00042670830161\\
8	1.00042670830158\\
16	1.00042670830158\\
32	1.00042670830158\\
64	1.00042670830158\\
128	1.00042670830158\\
256	1.00042670830158\\
512	1.00042670830158\\
1024	1.00042670830158\\
2048	1.00042670830158\\
4096	1.00042670830158\\
8192	1.00042670830158\\
16384	1.00042670830158\\
};
\addlegendentry{$p-p_0 \text{ in } L^\infty(L^2)$}

\addplot [color=mycolor4, mark=triangle, line width = 1.2]
  table[row sep=crcr]{%
1	1.00001809272496\\
2	1.00096857534677\\
4	1.0019241157644\\
8	1.00653181723556\\
16	1.01510661407494\\
32	1.0257630491871\\
64	1.03001277569842\\
128	1.02859906658891\\
256	1.02590711008459\\
512	1.02347960957052\\
1024	1.02151571515296\\
2048	1.01994070002055\\
4096	1.01902323137766\\
8192	1.01815087818523\\
16384	1.01748833017674\\
};
\addlegendentry{$p-p_0 \text{ in } L^2(H^1)$}

\addplot [color=mycolor2, mark=diamond, line width = 1.2, densely dashed, mark options={solid}]
  table[row sep=crcr]{%
1	1.46396151364356\\
2	1.42596141801432\\
4	1.35296304281532\\
8	1.23778823298487\\
16	1.10332796245079\\
32	1.04913797414562\\
64	1.03123982103674\\
128	1.02636615491439\\
256	1.02474272571448\\
512	1.02429875006148\\
1024	1.02418665699256\\
2048	1.02415868100518\\
4096	1.0241522084269\\
8192	1.02415074328109\\
16384	1.02415039648877\\
};
\addlegendentry{$\sqrt \eps (m-m_0) \text{ in } L^\infty(L^2)$}

\addplot [color=mycolor5, mark=square, line width = 1.2, dashed, mark options={solid}]
  table[row sep=crcr]{%
1	1.00057931187375\\
2	1.00033743801173\\
4	0.999616603522896\\
8	0.996698039170977\\
16	0.990170745158128\\
32	0.977327857089737\\
64	0.963859759686008\\
128	0.960091434504141\\
256	0.957366736553151\\
512	0.9566579915376\\
1024	0.956488355754195\\
2048	0.956448778456474\\
4096	0.956439551084165\\
8192	0.956437398645655\\
16384	0.956436889703709\\
};
\addlegendentry{$m-m_0 \text{ in } L^2(L^2)$}

\addplot [color=mycolor7, mark=o, line width = 1.2, dotted, mark options={solid}]
  table[row sep=crcr]{%
1	1.00057931187375\\
2	1.00033743801168\\
4	0.992373305368376\\
8	0.957084050867703\\
16	0.886018993941131\\
32	0.841957849980268\\
64	0.804873123296446\\
128	0.78839684414808\\
256	0.77732827910548\\
512	0.774539311357401\\
1024	0.77655609057319\\
2048	0.780197878433649\\
4096	0.782657633939377\\
8192	0.782682254196111\\
16384	0.782693412212063\\
};
\addlegendentry{$\lambda-\lambda_0 \text{ in } L^2(\R^2)$}
\end{axis}

\end{tikzpicture}%
	\caption{Estimate of the approximation order~$\alpha$ corresponding to~Theorems~\ref{thm_p0m0} and~\ref{thm_p0m0_2} with consistent initial values given in~\eqref{eqn:initialData:c}, i.e., with $m(0)=m_0(0)$. }
	\label{fig_p0m0_consistent}
\end{figure}
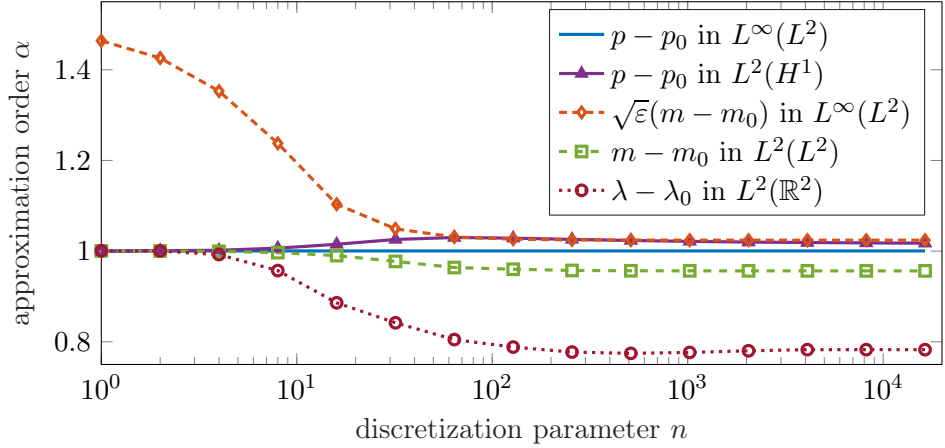

In Section~\ref{sect:expansion:first} we also proved that the convergence orders improve by half an order if $m(0)=m_0(0)$ or, equivalently, $\partial_x p(0) = - m(0)$ is satisfied. This is numerically confirmed in Figure~\ref{fig_p0m0_consistent}, 
where the associated initial values are given by 
\begin{align}
\label{eqn:initialData:c}
p(x,0)= \sum_{k=1}^\infty \frac{\sin(\pi k x)}{ k^{2.55}} \in H^2_0(0,1)
\quad\text{and}\quad 
m(x,0) = -\pi \sum_{k=1}^\infty \frac{\cos(\pi k x)}{k^{1.55}} \in H^1(0,1).
\end{align} 
At this point, we would like to emphasize that this also improves the rate of the Lagrange multiplier~$\lambda$ to~$0.78$ for larger $n$. With similar arguments as made for the first example, we would expect a rate~$0.75$. 

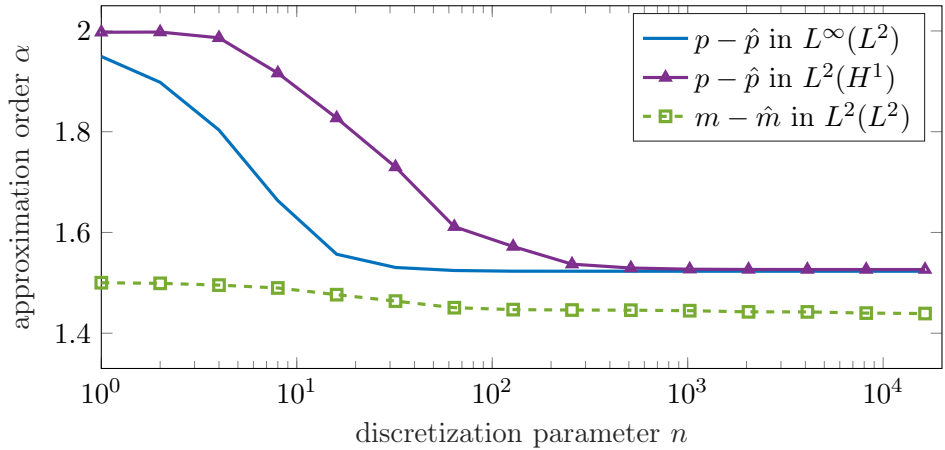
\begin{figure}
%
%
\definecolor{mycolor1}{rgb}{0.00000,0.44700,0.74100}%
\definecolor{mycolor2}{rgb}{0.85000,0.32500,0.09800}%
\definecolor{mycolor3}{rgb}{0.92900,0.69400,0.12500}%
\definecolor{mycolor4}{rgb}{0.49400,0.18400,0.55600}%
\definecolor{mycolor5}{rgb}{0.46600,0.67400,0.18800}%
\definecolor{mycolor6}{rgb}{0.30100,0.74500,0.93300}%
\definecolor{mycolor7}{rgb}{0.63500,0.07800,0.18400}%
\begin{tikzpicture}

\begin{axis}[%
width=4.4in,
height=1.9in,
at={(2.502in,1.239in)},
scale only axis,
xmode=log,
xmin=1,
xmax=20000,
xminorticks=true,
xlabel style={font=\color{white!15!black}},
xlabel={discretization parameter $n$},
ymin=1.33,
ymax=2.05,
ylabel style={font=\color{white!15!black}},
ylabel={approximation order $\alpha$},
axis background/.style={fill=white},
title style={font=\bfseries},
legend style={legend cell align=left, align=left, draw=white!15!black}
]
\addplot [color=mycolor1, line width = 1.2]
  table[row sep=crcr]{%
1	1.949337974485348\\
2	1.898060020610476\\
4	1.803311414322690\\
8	1.663174009508613\\
16	1.556772229276295\\
32	1.530463900757669\\
64	1.524471198692384\\
128	1.523094029037946\\
256	1.523048705633808\\
512	1.523046074073967\\
1024	1.523045905911523\\
2048	1.523045896147684\\
4096	1.523045895619877\\
8192	1.523045895588222\\
16384	1.523045895586300\\
};
\addlegendentry{$p-\hat p \text{ in } L^\infty(L^2)$}

\addplot [color=mycolor4, mark=triangle, line width = 1.2]
  table[row sep=crcr]{%
1	1.997415516923470\\
2	1.997807834785766\\
4	1.986344464545082\\
8	1.916302123865573\\
16	1.827095556114595\\
32	1.729955224464728\\
64	1.611436207885850\\
128	1.572076069395302\\
256	1.537085046323380\\
512	1.529430348185340\\
1024	1.527038387860456\\
2048	1.526529664191256\\
4096	1.526410894780790\\
8192	1.526383181922729\\
16384	1.526376876441625\\
};
\addlegendentry{$p-\hat p \text{ in } L^2(H^1)$}

\addplot [color=mycolor5, mark=square, line width = 1.2, dashed, mark options={solid}]
  table[row sep=crcr]{%
1	1.500212332387225\\
2	1.498985200250401\\
4	1.495494812002772\\
8	1.489823919263869\\
16	1.476484120862075\\
32	1.463740472891805\\
64	1.450655845633195\\
128	1.446841972827847\\
256	1.446032782479371\\
512	1.445718023846636\\
1024	1.444674415328883\\
2048	1.442439170609907\\
4096	1.442274763411707\\
8192	1.439957452640277\\
16384	1.439044443387655\\
};
\addlegendentry{$m-\hat m \text{ in } L^2(L^2)$}

\end{axis}

\end{tikzpicture}%
	\caption{Estimate of the approximation order~$\alpha$ corresponding to~Theorem~\ref{thm_p1m1} with the initial values from~\eqref{eqn:initialData:d}.}
	\label{fig_p1m1}
\end{figure}

Finally, to verity the results of Theorem~\ref{thm_p1m1}, i.e., considering $\hat p$ and $\hat m$, we set as initial data 
\begin{align}
\label{eqn:initialData:d}
p(x,0)= \sum_{k=1}^\infty \frac{\sin(\pi k x)}{ k^{3.55}} \in H^3_0(0,1)
\quad\text{and}\quad 
m(x,0) = -\pi \sum_{k=1}^\infty \frac{\cos(\pi k x)}{k^{2.55}} \in H^2(0,1). 
\end{align}
These values are obviously consistent, i.e., $\partial_x p(0) = - m(0)$, and  satisfy the smoothness requirements of Theorem~\ref{thm_p1m1}. The resulting approximations of the rates are displayed in Figure~\ref{fig_p1m1}. 
Here, the approximated order for the pressure fades from two to around~$1.525$, which is the expected rate. As before, the mass flux shows a slightly worse result than proven in Theorem~\ref{thm_p1m1}, presumably due to round off errors. A rigorous calculation proves the predicted order of $1.5$ in this case. 

In summary, apart from small discrepancies due to round-off errors, all $\eps$-rates for the pressure and the Lagrange multiplier show a decreasing asymptotic behavior in the discretization parameter $n$. The approximate orders for the mass flux $m$, on the other hand, are more or less independent of~$n$ as predicted in Section~\ref{sect:expansion}.

%

%
%
\section{Conclusion}\label{sect:conclusion}
In this paper, we have considered linear PDAEs of hyperbolic type with a small parameter~$\eps>0$, which turn parabolic in the limit case, i.e., for $\eps=0$. Depending on the consistency of the initial data and the regularity of the right-hand sides, we have shown first- and second-order estimates of the corresponding expansion in terms of~$\eps$. In a number of numerical experiments, we have validated these results and compared them with the finite-dimensional setting. 
The presented expansion may be used for the construction of novel numerical methods. For this, the approach needs to be combined with integration schemes for PDAEs of parabolic type such as splitting schemes~\cite{AltO17}, Runge-Kutta methods~\cite{AltZ18,Zim21}, discontinuous Galerkin methods~\cite{VouR18}, or exponential integrators~\cite{AltZ20,Zim21}. 
%
%
%
%
\bibliographystyle{alpha} 
\newcommand{\etalchar}[1]{$^{#1}$}

%
%
%
\appendix
\section{Proofs}\label{appendix}
Within this appendix, we collect several proofs and start with showing that the operator~$A_\gamma$ generates a $C_0$-semigroup. 
%
\begin{proof}[Proof of Lemma~\ref{lem_unbounded_A}]
Without loss of generality, we assume that $\A$ is non-negative. Otherwise, we consider $A_{\gamma}-\gamma C_{\A}\id_{\cHker\times\M}$ and use \cite[p.~12]{Paz83}. This then updates the operator $\D$ in the $(2,2)$-component of~$A_{\gamma}$ to $\D + \gamma^2 C_{\A}\id_{\M}$, which is still elliptic on $\M$ for every $\gamma>0$. 

Since the term~$\A q_{\ker}$ is an element of~$\cH^\ast \subseteq \cHker^\ast \cong \cHker$ for every $q_{\ker} \in \Pker$, the operator $A_\gamma$ is bounded on its domain $D(A_\gamma)$. We show that $A_{\gamma}$ is a densely defined, closed, and dissipative operator with a dissipative adjoint~$A_{\gamma}^*$. The statement then follows by~\cite[Ch.~1.4, Cor.~4.4]{Paz83}. 
By Lemma~\ref{lem_op_B_station} the operator $A_\gamma$ is linear, bounded, maps from $\Pker \times \calM$ to $\Pker^\ast \times \calM^\ast$, and has a bounded inverse. 
In particular, it holds that $A_{\gamma}^{-1} (g,f) \in D(A_\gamma)$ for all $(g,f)\in \cHker\times \calM \hookrightarrow \Pker^\ast \times \calM^\ast$. 
This proves the closeness by a simple calculation. 
Furthermore, the operator~$A_{\gamma}$ is dissipative by Definition~4.1 in~\cite[Ch.~1.4]{Paz83}. To see this, consider a given $(p,m)\in D(A_\gamma) \hookrightarrow \Pker\times \calM$ and choose the same element as test function under the embedding $\Pker\times \calM \hookrightarrow \cHker\times \calM$. The adjoint $A_{\gamma}^\ast$ is dissipative as well, since the adjoint operators~$\A^\ast\colon \P \to \P^\ast$ and~$\D^\ast\colon \calM \to \calM^\ast$ have the same properties as $\A$ and $\D$, respectively. 

It remains to show that~$A_\gamma$ is densely defined. Since~$D(A_\gamma)$ is independent of $\gamma$, we may fix~$\gamma =1$ for the remainder of the proof. Let $(h_{\ker},m)\in \calH_{\ker}\times \M$ be arbitrary. By the embeddings given by the Gelfand triple $\Pker, \calH_{\ker}, \Pker^*$, there exist for every~$\delta >0$ elements $p_{\ker} \in \Pker$ and $\g^{\prime} \in \calH_{\ker}$ with $\|p_{\ker} - h\|_{\calH} < \delta$ and $\|\g^{\, \prime}  - (\A p_{\ker}  - \K^* m) \|_{\Pker^*} < \delta$. 
Let $(p_{\ker}^{\, \prime}, m^{\prime}) \in \Pker \times \M$ be the unique solution of 
\begin{equation*}
\begin{alignedat}{5}
\A  p_{\ker}^{\, \prime} &\ -\ &\K^* m^{\prime} &=&&\ \g^\prime &&\qquad \text{in }\P^*_{\ker},\\
\K p_{\ker}^{\, \prime} &\ +\ &\D  m^{\prime} &=&&\ \K p_{\ker} + \D  m  &&\qquad \text{in }\M^*.
\end{alignedat}
\end{equation*}
By construction, we then have $(p_{\ker}^{\, \prime}, m^{\prime}) \in D(A_1)=D(A_\gamma)$. We finally choose $(p_{\ker}^{\, \prime}, m^{\prime})$ as approximation of $(h,m)$ and conclude with the boundedness of~$-A^{-1}_{1}$ that 
\begin{align*}
\| h-p_{\ker}^{\, \prime}\|_{\calH} + \| m - m^{\prime}\|_{\M} 
&\lesssim \| h-p_{\ker}\|_{\calH} + \|p_{\ker} - p_{\ker}^{\, \prime}\|_{\P} + \| m -m^{\prime}\|_{\M} \\
&\lesssim \| h-p_{\ker}\|_{\calH} + \| \g^{\prime}  - (\A p_{\ker}  - \K^* m) \|_{\Pker^*} 
< 2\, \delta. \qedhere
\end{align*}
\end{proof}
We now turn to the existence of mild and classical solutions of the considered PDAE model. 
%
%
\begin{proof}[Proof of Proposition~\ref{prop_mild_sol_B}]
Recall the equations for $(\pt, \mt, \lt)$ including the constraint~$\B \pt = 0$. This allows to reduce the ansatz and test space accordingly, leading to the equivalent system 
\begin{subequations}
	\label{eqn_inproof_ph}
	\begin{alignat}{5}
	\dot \pt &\ +\ & \A\pt &\ -\ & \K^* \mt &= \g_2 + C_{\A} \pb  - \dot \pb - \A\B^- \h - \B^- \dot \h  && \qquad \text{in } \Pker^*, \label{eqn_inproof_ph_a}\\
	\eps\, \dot \mt &\ +\ & \K \pt &\ +\ & \D \mt &= \f - \K \B^- \h - \eps\, \dot \mb  && \qquad \text{in } \M^*. \label{eqn_inproof_ph_b}
	\end{alignat}
\end{subequations}
At this point, we would like to emphasize that an equation stated in~$\Pker^*$ means that we only consider test functions in~$\Pker$. 
With the state $x := [\tfrac{1}{\sqrt{\eps}} \pt, \mt]^T$, equation~\eqref{eqn_inproof_ph} becomes the abstract Cauchy problem 
\begin{subequations}
	\label{eqn_cauchy_on_ker}
	\begin{align}
	\dot{x}&= \tfrac{1}{\sqrt{\eps}} A_\gamma x + F =
	\tfrac{1}{\sqrt{\eps}} A_\gamma x + \begin{bmatrix}
	\tfrac{1}{\sqrt{\eps}} (\g_2 + C_{\A} \pb - \dot \pb - \A\B^-\h - \B^- \dot \h)\\
	\tfrac{1}{\eps} (\f - \K \B^- h) - \dot \mb
	\end{bmatrix} \label{eqn_cauchy_on_ker_a}
	\end{align}
	with initial condition
	\begin{align}
	x(0) &= \big[\tfrac{1}{\sqrt{\eps}} \pt(0),\, \mt(0) \big]^T.  \label{eqn_cauchy_on_ker_b}
	\end{align}
\end{subequations}
Here, $A_\gamma$ equals the operator from Lemma~\ref{lem_unbounded_A} and $\gamma = \sqrt{\eps}$. 
Since the right-hand side satisfies $F\in L^2(0,T; \cH^* \times \M^*)$ and $x(0) \in \cHker \times \M$,  the Cauchy problem has a unique mild solution $x \in C([0,T],\cHker \times \M)$. 
Thus, $p = \pt + \pb + \B^-h \in C([0,T], \cH)$ has a derivative in $L^2(0,T;\Pker^*)$ by~\eqref{eqn_inproof_ph_a} and $m=\mt + \mb$ is an element of $C([0,T], \M)$. Finally, $\lambda$ can be constructed as in the proof of~\cite[Th~3.3]{EmmM13}. 
\end{proof}
%
%
%
\begin{proof}[Proof of Proposition~\ref{prop_clas_sol_B}]
Consider the auxiliary problem~\eqref{eq:op:stationary} with right-hand side~$\g_1$, leading to~$(\pb,\mb,\lb)  \in H^2(0,T;\Pker) \times H^2(0,T;\M) \times H^2(0,T;\Q)$, cf.~Lemma~\ref{lem_op_B_station}. 
Following the proof of Proposition~\ref{prop_mild_sol_B}, we notice that the right-hand side of the Cauchy problem~\eqref{eqn_cauchy_on_ker} is an element of~$H^1(0,T;\cH^\ast\times \M^\ast)$ by the more regular right-hand sides. 
For the initial data from~\eqref{eqn_cauchy_on_ker_b}, we note that~$x(0) \in D(A_{\sqrt\eps})$, since $\pt(0) \in \Pker$ and $\mt(0)$ satisfies 
\[
\K^* \mt(0)
= \K^* m(0) - \K^* \mb(0)
= \K^* m(0) + \g_1(0) - \A \pb(0)
= \mstar - \A \pb(0)
\quad\text{in } \cHker^\ast 
\]
due to \eqref{eq:op:stationary:a}. This means that~$\K^* \mt(0)$ has a representation in~$\cHker$. 
The claimed solution spaces of~$p$ and~$m$ follow by~\cite[Ch.~4, Cor.~2.10]{Paz83}.  
Since the operator~$\B$ satisfies an inf-sup condition by assumption, there exists a unique (and continuous) multiplier~$\lambda$, which satisfies~\eqref{eq:PDAE:eps:a}, cf.~\cite[Lem.~III.4.2]{Bra07}. 
\end{proof}
\end{document}